\documentclass[11pt,reqno]{amsart}
\usepackage{amsfonts}
\usepackage{mathrsfs}
\usepackage{graphicx} 
\usepackage{epsfig}
\usepackage{amsmath, amsthm}
\usepackage{amssymb}

\theoremstyle{plain}
\newtheorem{thm}{Theorem}[section]
\newtheorem{prop}[thm]{Proposition}
\newtheorem{lem}[thm]{Lemma}
\newtheorem{cor}[thm]{Corollary}

\theoremstyle{definition}

\newtheorem{rem}[thm]{Remark}
\newtheorem{defn}[thm]{Definition}
\newtheorem{eg}[thm]{Example}
\newtheorem{subtitle}[thm]{}
\newtheorem{ex}{Exercise}[section]
\numberwithin{equation}{section}

\def\a{\alpha}
\def\b{\beta}
\def\d{\delta}

\def\g{\gamma}

\def\K{\nabla}
\def\l{\lambda}

\def\n{\vert\/}
\def\o{\theta}

\def\li{\langle}
\def\ri{\rangle}
\def\n{|\/ }
\def\tr{{\rm tr}}
\def\bs{\bigskip}
\def\ms{\medskip}
\def\ss{\smallskip}

\def\di{$\diamond$}
\def\ni{\noindent}
\def\ti{\tilde}
\def\p{\partial}
\def\Re{{\rm Re\/}}
\def\Im{{\rm Im\/}}
\def\I{{\rm I\/}}

\def\diag{{\rm diag}}
\def\ad{{\rm ad}}
\def\Ad{{\rm Ad}}
\def\Iso{{\rm Iso}}
\def\Gr{{\rm Gr}}

\def\C{\mathbb{C}}

\def\H{\mathbb{H}}

\def\R{\mathbb{R} }

\newcommand{\beq}{\begin{equation}}
\newcommand{\eeq}{\end{equation}}
\newcommand{\beg}{\begin{eg}}
\newcommand{\eeg}{\end{eg}}
\newcommand{\bthm}{\begin{thm}}
\newcommand{\ethm}{\end{thm}}
\newcommand{\bprop}{\begin{prop}}
\newcommand{\eprop}{\end{prop}}
\newcommand{\bcor}{\begin{cor}}
\newcommand{\ecor}{\end{cor}}
\newcommand{\blem}{\begin{lem}}
\newcommand{\elem}{\end{lem}}
\newcommand{\bca}{\begin{cases}}
\newcommand{\eca}{\end{cases}}
\newcommand{\brem}{\begin{rem}}
\newcommand{\erem}{\end{rem}}
\newcommand{\bpm}{\begin{pmatrix}}
\newcommand{\epm}{\end{pmatrix}}
\newcommand{\bbm}{\begin{bmatrix}}
\newcommand{\ebm}{\end{bmatrix}}
\newcommand{\bvm}{\begin{vmatrix}}
\newcommand{\evm}{\end{vmatrix}}
\newcommand{\bdefn}{\begin{defn}}
\newcommand{\edefn}{\end{defn}}
\newcommand{\bsub}{\begin{subtitle}}
\newcommand{\esub}{\end{subtitle}}
\newcommand{\bex}{\begin{ex}}
\newcommand{\eex}{\end{ex}}
\newcommand{\ben}{\begin{enumerate}}
\newcommand{\een}{\end{enumerate}}

\date{} 

\def\calA{\mathcal{A}}

\def\calC{\mathcal{C}}

\def\calG{\mathcal{G}}
\def\calI{\mathcal{I}}

\def\calK{\mathcal{K}}

\def\calM{\mathcal{M}}
\def\calN{\mathcal{N}}
\def\calO{\mathcal{O}}
\def\calP{\mathcal{P}}
\def\R{\mathbb{R}}
\def\calS{\mathcal{S}}

\def\calU{\mathcal{U}}

\def\R{\mathbb{R}}

\def\C{\mathbb{C}}

\def\fr{\mathfrak{r}}
\def\fR{\mathfrak{R}}
\def\fy{\mathfrak{y}}

\def\rg{\rm g}

\def\diag{{\rm diag \/ }}
\def\det{{\rm det \/ }}

\def\rc{\R^n \backslash \{0\}}
\def\bh{\backslash}

\def\rd{{\rm \/ d\/}}
\def\bu{\bullet}

\def\di{\diamond}

\def\an1{A^{(1)}_{n-1}}

\def\Diff{{\rm Diff\,}}

\def\r2{\R^2\backslash \{0\}}

\def\rb{\bf b}
\def\rn{\bf n}

\def\n{\, | \,}

\def\di{\diamond}

\begin{document}

\title[Dispersive Geometric Curve flows]
{Dispersive Geometric Curve Flows} 

\author{Chuu-Lian Terng$^\dag$}\thanks{$^\dag$Research supported
in  part by NSF Grant DMS-1109342}
\address{Department of Mathematics\\
University of California at Irvine, Irvine, CA 92697-3875.  Email: cterng@math.uci.edu}


\maketitle

\begin{abstract}
 The Hodge star mean curvature flow on a 3-dimension Riemannian or pseudo-Riemannian manifold, the geometric Airy flow on a Riemannian manifold, the Schr\"dingier flow on Hermitian manifolds, and the shape operator curve flow on submanifolds are natural  non-linear dispersive curve flows in geometric analysis.  A curve flow is integrable if the evolution equation of the local differential invariants of a solution of the curve flow is a soliton equation. For example, the Hodge star mean curvature flow on $\R^3$ and on $\R^{2,1}$, the geometric Airy flow on $\R^n$, the Schr\"dingier flow on compact Hermitian symmetric spaces, and the shape operator curve flow on an Adjoint orbit of a compact Lie group are integrable. In this paper, we give a survey of these results, describe a systematic method to construct integrable curve flows from  Lax pairs of soliton equations, and discuss the Hamiltonian aspect and the Cauchy problem of these curve flows.

\end{abstract}

\section{Introduction}

Three of the simplest linear dispersive equations in one space and one time are:
\ben
\item[$\bu$] {\it The  linear Schr\"dingier equation\/} for $q:\R^2\to \C^n$,
 $$q_t= i q_{xx}.$$ 
 \item[$\bu$] {\it The Airy equation\/} for $q:\R^2\to \R$,
$$q_t= q_{xxx}.$$ 
\item[$\bu$] {\it The first order symmetric linear hyperbolic system\/} for  $q:\R^2\to \R^n$, 
$$q_t= A q_x,$$ where $A$ is a constant $n\times n$ real symmetric matrix.
\een
These are Hamiltonian partial differential equations, and the Cauchy problem can be solved by Fourier transforms.  There are also many interesting completely integrable non-linear dispersive Hamiltonian partial differential equations (soliton equations) that arise in the study of non-linear waves. For example, 
\ben
\item[$\bu$]  the {\it focusing non-linear  Schr\"dinger equation\/} (NLS), 
$$q_t=\frac{i}{2}(q_{xx}+ 2|q|^2 q),$$ for $q:\R^2\to \C$, 
\item[$\bu$]  the {\it sine-Gordon equation\/} (SGE) for $q:\R^2\to \R$,
$$q_{tt}-q_{xx}=\sin q,$$
\item[$\bu$]  the {\it KdV equation\/}
$$q_t=\frac{1}{4}( q_{xxx}+6qq_x),$$
\item[$\bu$]  the {\it $n$ wave equation\/} for $u=(u_{ij}):\R^2\to u(n)$ with $u_{ii}=0$ for $1\leq i\leq n$,
$$(u_{ij})_t= \frac{b_i-b_j}{a_i-a_j}(u_{ij})_x +\sum_{k\not= i, j}\left(\frac{b_k-b_j}{a_k-a_j}-\frac{b_i-b_k}{a_i-a_k}\right) u_{ik} u_{kj}, \quad i\not=j$$
where $a_i, b_i$ are complex constants for $1\leq i\leq n$ and $a_1, \ldots, a_n$ are distinct.  
\een

Soliton equations have many remarkable properties including:
\ben
\item[$\bu$] There are infinitely many families of explicit soliton solutions.
\item[$\bu$] There exist B\"acklund transformations that generate new solutions from a given solution.
\item[$\bu$] There exists a {\it Lax pair\/}, i.e., a one parameter family of $\calG$-valued connection one forms $\o(x,t,\l)=A(u, \l) \rd x+ B(u,\l) \rd t$ on the $(x,t)$ plane defined by $u$, $x$ derivatives of $u$ and a parameter $\l$ such that $\o(\cdot, \cdot, \l)$ is flat for all $\l\in \C$ if and only if $u$ is a solution of the soliton equation.
\item[$\bu$] There are interesting actions of infinite dimension groups on the space of solutions.
\item[$\bu$] There exists a bi-Hamiltonian structure, i.e., two compatible Poisson structures such that the soliton equation are Hamiltonian with respect to both structures.
\item[$\bu$] There are infinitely many commuting conservation laws. 
\item[$\bu$] If $u$ is a solution of the soliton equation then the scattering of the linear operator $\rd_x + A(u(\cdot, t),\l)=0$ (the $x$ part of the Lax pair) evolves linearly in $t$. If the inverse scattering transform of this linear operator exists for a class $\calO$ of initial data then the Cauchy problem with initial data in $\calO$ of the soliton equation can be solved.
\een

One of the standard methods for constructing soliton equations is to use splittings of infinite dimensional Lie algebras, i.e., the Adler-Kostant-Symes construction.  The properties of soliton equations given above can be derived in a unified and systematic way from these splittings (cf. \cite{DS84}, \cite{TerUhl98}). 

Soliton equations also arise naturally in differential geometry (cf. \cite{Ter03}). For example, the SGE is the equation for surfaces in $\R^3$ with Gaussian curvature $-1$, and  the reduced $n$-wave equation is related to the Gauss-Codazzi equations for flat Lagrangian submanifolds in $\C^n$ with flat and non-degenerate normal bundle (cf. \cite{TerWan08}).  

The linear dispersive equations mentioned above are flows on the space $C^\infty(\R, V)$ of smooth maps from $\R$ to $V$, where $V$ is $\R^n$ or $\C^n$.  It is natural in geometric analysis to construct non-linear analogues of linear dispersive equations by replacing $V$ by a Riemannian, pseudo-Riemannian, or a Hermitian manifold, and replacing $\p_x$ by $\K_{e_0}$ or $\K_{\g_x}$, where $\K$ is the metric connection and $e_0$ is the unit tangent vector along the curve.
 Next we give some natural analogues of linear dispersive equations in differential geometry. 

\ben
\item[$\bu$] The  {\it Schr\"odingier flow\/} on Hermitian manifolds

Let  $(M, \rg, J)$ be a Hermitian manifold, where $J$ is the complex structure and $\rg$ is a Hermitian metric. The Schr\"odinger flow on $M$  is the following natural analogue of the linear Schr\"dingier equation,
\beq\label{be}
\g_t= J_\g(\K_{\g_x}\g_x),
\eeq
where $\K$ is the Levi-Civita connection of $\rg$. 

\item[$\bu$] {\it The star mean curvature flow ($\ast$-MCF) on a 3-manifold} 

Let $\rg$ be a Riemannian or Lorentzian metric on a $3$-dimension manifold $N$. 
The $\ast$-MCF on $(N,\rg)$ is the following curve flow on the space of immersed curves in $N$,
$$\g_t= \ast_\g(H(\g(\cdot, t))),$$
where $\ast_{\g(x)}$ is the Hodge star operator on the normal plane $\nu(\g)_x$ and $H(\g(\cdot, t))$ is the mean curvature vector for $\g(\cdot, t)$.

\item[$\bu$] {\it The geometric Airy flow on a Riemannian manifold}

Let $(M,\rg)$ be a Riemannian manifold, $\K$  the Levi-Civita connection of $\rg$, and $\K^\perp$ the induced normal connection on a curve $\g$ defined by 
$\K^\perp_{e_0} \xi=( \K_{e_0}\xi)^\perp$, the projection of $\K_{e_0}\xi$ onto $\nu(\g)$, where $\xi$ is a normal vector field $\xi$ along $\g$ and $e_0$ is the unit tangent.  The geometric Airy flow on $(M,\rg)$ is
$$\g_t= \K^\perp_{e_0}H(\g(\cdot, t)).$$

\item[$\bu$] {\it The shape operator curve flow on a submanifold}

Let $M$ be a submanifold of a Riemannian manifold $(N,\rg)$, and $\eta$ a normal field on $M$.  The shape operator $A_{\eta(x)}:TM_x\to TM_x$ is the self-adjoint operator defined by 
$$A_\eta(v)= -(\K_v \eta)^T,$$
the tangential component of $\K_v\eta$, where $\K$ is the Levi-Civita connection of $\rg$.
{\it The shape operator curve flow on $M$} is 
$$\g_t= A_{\eta}(\g_x).$$  This is a natural non-linear symmetric first order hyperbolic system on the submanifold $M$. 
\een

Given a geometric curve flow on a homogeneous $G$-space $M$, if we can construct a ``good'' moving frame along curves on $M$ such that the evolution equation of the differential invariants of solutions of a geometric curve flow defined by these moving frames is a soliton equation, then we can use techniques from soliton theory to study this curve flow.  We call such a curve flow an {\it integrable curve flow\/}. As we will see next, many of the geometric dispersive curve flows are known to be integrable.

It can be checked that the $\ast$-MCF preserves the arc-length parameter. If $\g$ is parametrized by its arc-length, then the mean curvature vector of $\g$ is $H(\g)= \K_{\g_x}\g_x$.  So the $\ast$-MCF on $\R^3$ parametrized by the arc-length is 
$$\g_x= \ast_\g H(\g(\cdot, t))= \g_x\times H(\g(\cdot, t))= \g_x\times \g_{xx}.$$  So the $\ast$-MCF on $\R^3$ is the {\it vortex filament equation\/} (VFE),
\beq\label{dx}
\g_t= \g_x\times \g_{xx}.
\eeq
Hasimoto proved in \cite{Has72} that if $\g$ is a solution of the VFE parametrized by arc-length then there exists a smooth $\o:\R\to \R$ such that 
$$q(x,t)= \frac{1}{2}e^{i\o(t)} k(x,t) e^{-\int_0^x \tau(s, t)\rd s}$$
is a solution of the NLS, where $k(\cdot, t)$ and $\tau(\cdot, t)$ are the curvature and torsion along $\g(\cdot, t)$. 
Note that
\ben
\item[(i)] the VFE is invariant under the group $\fR(3)$ of rigid motions of $\R^3$, i.e., if $\g$ is a solution of the VFE, then so is $C\circ \g$ for any $C\in \fR(3)$,
\item[(ii)] if $q$ is a solution of the NLS and $c\in \C$ with $|c|=1$, then $cq$ is also a solution of the NLS, i.e., the group $S^1$ acts on the space of solutions of the NLS,
\item[(iii)] if $C\in\fR(3)$, then $\g$ and $C\g$ have the same curvature and torsion.
\een
So the {\it Hasimoto's transform\/} in fact maps the $\fR(3)$-orbit of a solution of the VFE to the $S^1$-orbit of a solution of the NLS. We use $\Psi$ to denote the Hasimoto transform on these orbit spaces, i.e.,
\begin{align*}
&\{\g:\R^2\to \R^3\n \g \, {\rm\, is\, a\, solution\, of\, the\,VFE,\, \,} ||\g_x||=1\}/\fR(3)\\
&\hskip 1.5in \downarrow\, \Psi\\
&\qquad\{q:\R^2\to \C\n q \, {\rm is \, a\, solution\, of\, the\, NLS\,}\}/S^1
\end{align*}
 The inverse of $\Psi$ is given by the Pohlmeyer-Sym construction (cf. \cite{Poh76}, \cite{Sym85}).  
 
 Here is an outline of the paper:
 \ben
 \item[(a)]  An orthonormal frame $(e_0, e_1, e_2)$ along a curve $\g$ in $\R^3$ is a {\it p-frame\/} if $e_0$ is the unit tangent and $(e_1, e_2)$ is a parallel orthonormal frame of the normal bundle $\nu(\g)$ with respect to the induced normal connection. The differential invariants $k_1, k_2$ of $\g$ defined by a p-frame are the {\it principal curvatures\/} along $e_1, e_2$.  Hasimoto's result can be stated as follows: If $\g:\R^2\to \R^3$ is a solution of the VFE parametrized by the arc-length, then there exists $g:\R^2\to SO(3)$ such that $g(\cdot, t)$ is a p-frame along $\g(\cdot, t)$ for all $t\in \R$ and $q=\frac{1}{2}(k_1+ ik_2)$ is a solution of the NLS, where $k_1, k_2$ are the principal curvatures defined by the p-frame $g$.   Although p-frames along $\g$ are not unique and principal curvatures of $\g$ depend on the choice of p-frames, the formula for the Hasimoto transform $\Psi$ takes a simpler form in p-frames. It is also easier to use p-frames to translate the results in soliton theory for the NLS to the VFE. The VFE is one of the most well-known integrable curve flow, and there are many works on the relation between the VFE and NLS (cf. \cite{LanPer91}, \cite{DoSa94}, \cite{SaYa98}, \cite{CaIv05}, \cite{CaIv07}). We use the VFE as a model example to explain how to use the Hasimoto transform $\Psi$ to study the bi-Hamiltonian structure, higher order commuting curve flows, and B\"acklund transformations for the VFE.  We also show that the geometric Airy flow on $\R^3$ preserves the total arc-length. Hence we can reparametrize the curves by the arc-length parameter, and the resulting normalized geometric Airy flow commutes with the VFE. 
\item[(b)] We use p-frames to construct analogues of the Hasimoto transform $\Psi$ between the time-like $\ast$-MCF on $\R^{2,1}$ and the defocusing NLS, and between the space-like $\ast$-MCF on $\R^{2,1}$ and the  $2\times 2$ AKNS equation
\beq\label{av}
 \bca q_t= -\frac{1}{2}(q_{xx}- 2q^2r),\\ r_t= \frac{1}{2}(r_{xx} - 2qr^2).\eca
 \eeq 
 The bi-Hamiltonian structure, higher order commuting curve flows, and B\"acklund transformations for the $\ast$-MCF on $\R^{2,1}$ are constructed.
\item[(c)]
Analogues of the Hasimoto transforms were also constructed between the Schr\"dingier flow on $\Gr(k,C^n)$ and the matrix NLS in \cite{TU06}, between  central affine curve flows on $\rc$ and the Gelfand-Dickey hierarchy in \cite{UP95} and \cite{CIM09} for $n=2$, in \cite{CIM13} for $n=3$, and in \cite{TW14b} for $n\geq 3$, and between the shape operator curve flow on Adjoint orbits of $U(n)$ on $u(n)$ and the $n$-wave equation in \cite{Fer95} and \cite{TerTho01}.  We will give a brief survey of these results in the last section.
\item[(d)] These analogous Hasimoto transforms are constructed by choosing suitable moving frames so that equations for the differential invariants are soliton equations. A $p$-frame for a curve in $\R^3$ is only unique up to an action of $H= SO(2)$. This  motivates us to define a notion of $(H,V)$-moving frames for curves on homogeneous $G$-spaces, where $V$ is an affine subspace of $\calG$ where the differential invariants lies. All moving frames used in examples given in (a)-(c) are $(H,V)$-moving frames. Moreover, the notion of $(H,V)$-moving frame is equivalent to the definition of $H$-slice in the theory of transformation groups.
\item[(e)] Given a soliton equation, we give a general method for constructing geometric curve flows and suitable $(H,V)$-moving frames from the Lax pair so that their invariants give rise to solutions of the given soliton equation.  All integrable curve flows mentioned in this paper can be constructed by this method. 
\een

This paper is organized as follows: In section \ref{by}, we define the notion of $(H,V)$-moving frames for curves on a homogeneous space, give some examples, and explain the relation between $(H,V)$-moving frames and {\it $H$-slice\/} in the theory of transformation groups. We give a brief review of the AKNS $2\times 2$-hierarchy and its various restrictions in section \ref{bz}. The relation between the VFE and the NLS and the Cauchy problem for the VFE with initial data having rapidly decaying principal curvatures are discussed in section \ref{ca}. We solve the Cauchy problem with periodic initial data for the VFE in section \ref{du}. In section \ref{dv}, we construct B\"acklund transformations and give an algorithm to write down infinitely many families of explicit soliton solutions for the VFE. We explain the Hamiltonian aspects of the VFE including higher commuting curve flows and show that the normalized geometric Airy flow on $\R^3$ commutes with the VFE in section \ref{ef}.    In section \ref{dw}, we study the time-like $\ast$-MCF flow on $\R^{2,1}$. We study the space-like $\ast$-MCF on $\R^{2,1}$ and prove that the space-like geometric Airy flow on $\R^{2,1}$ with a suitable constraint gives a natural geometric interpretation of the KdV in section \ref{dwa}. We explain (e) in the last section.

\bs
\section{Moving frames along curves}\label{by}

We first give a review of local theory of curves in a Riemannian manifold, then define and give examples of $(H,V)$-moving frames for curves on a homogeneous space. In the end of the section we explain the relation between the notion of $(H,V)$-moving frames on a homogeneous space $G/K$ and $H$-slices for the gauge action of $C^\infty(\R, K)$ on the space of connections $\rd_x+ C^\infty(\R, \calG)$. 

Let $\g:\R\to M$ be a curve in a Riemannian manifold $(M^n,\rg)$, and $e_0$ the unit tangent along $\g$.  Let $g=(e_0, \ldots, e_{n-1})$ be an  orthonormal frame along $\g$ such that $e_0$ is the unit tangent to $\g$.  Then 
$$\K_{e_0} (e_0, \ldots, e_{n-1}):= (\K_{e_0} e_0, \ldots, \K_{e_0}e_{n-1}) =(e_0, \ldots, e_{n-1}) {\bf k}$$
for some ${\bf k}=(k_{ij})_{0\leq i, j\leq n-1}\in so(n)$, where $\K$ is the Levi-Civita connection for $\rg$.  In fact, $k_{ij}= \rg(\K_{e_0} e_j, e_i)$ for $0\leq i, j\leq n-1$, $k_i= k_{i0}=\rg(\K_{e_0} e_0, e_i)$ is the {\it principal curvature of $\g$\/} along $e_i$ for $1\leq i\leq n-1$, and
$$H(\g)= \K_{e_0}e_0=\sum_{i=1}^{n-1} k_i e_i$$
is  the {\it mean curvature vector field\/} along $\g$ and is independent of the choice of orthonormal frames. In particular, if $\g$ is parametrized by its arc-length, then 
$$H(\g)=\K_{\g_x} \g_x.$$
 
 Principal curvatures depend on the choice of normal frame. 
   If both $(e_0, \ldots, e_{n-1})$ and $(e_0, \ti e_1, \ldots, \ti e_{n-1})$ are orthonormal frames along $\g$ then there is a $h:\R\to SO(n-1)$ such that $(\ti e_1, \ldots, \ti e_{n-1})= (e_1, \ldots, e_{n-1}) h$. So the corresponding principal curvatures are related by
\beq\label{aw}
 (\ti k_1, \ldots, \ti k_{n-1})^t= h^{-1}(k_1, \ldots, k_{n-1})^t.
\eeq 
Note that an orthonormal frame $(e_1, \ldots, e_{n-1})$ of $\nu(\g)$ is parallel with respect to $\K^\perp$ if and only if $k_{ji}=\rg((e_i)_x, e_j)=0$ for all $1\leq i, j\leq n-1$. Moreover, two parallel normal frames are differed by a constant in $SO(n-1)$.   
Motivated by this example, we make the following definition.

\bdefn Let  $M=G\cdot p_0$ be a homogeneous $G$-space, $\calM$ a subset of $C^\infty(\R, M)$, $H$ a closed subgroup of  the isotropy subgroup  $G_{p_0}$ at $p_0$, and $V$ an affine subspace of $\calG$. We say that $\calM$ admits {\it $(H,V)$-moving frame\/} if 
\ben
\item given $\g\in \calM$, there exists  $g:\R\to G$ such that 
\ben
\item $\g(x)= g(x)\cdot p_0$,
\item $g^{-1}g_x\in C^\infty(\R, V)$.
\een
Moreover,  if $g_1\in C^\infty(\R, G)$ also satisfies (a)-(b), then there is a constant $h\in H$ such that $g_1= gh$.
\item if $g:\I\to G$ satisfying $g^{-1}g_x\in V$, then $\g= g\cdot p_0$ is in $\calM$. 
\een
We call such $g$ a {\it $(H,V)$-moving frame along $\g$\/} and $u=g^{-1}g_x$ the {\it differential invariants of $\g$\/} defined by $g$. 
\edefn

Note that if  both $g$ and $\ti g$ are $(H,V)$-frames along $\g$, then there is $c\in H$ such that $\ti g= gc$. So the corresponding differential invariants $\ti u= \ti g^{-1} g_x$ and $u=g^{-1}g_x$ are related by $\ti u= cuc^{-1}$.  

If $\calM$ admits $(H,V)$ moving frames, then by the existence and uniqueness of ordinary differential equations we have the following: 
\ben
\item[$\diamond$] $\calM$ is invariant under the action of $G$.
\item[$\diamond$] Differential invariants of $\g$ in $\calM$ determines $\g$ uniquely up to the action of $G$.
\een 
 
Next we give some examples. 

\beg[{\bf p-frame for curves in $\R^n$}]\ 

The group $\fR(n)$ of rigid motions of $\R^n$ can be embedded as the following subgroup of $GL(n+1)$:
$$\fR(n)=\left\{\bpm 1& 0\\ y & g\epm\,\bigg|\, y\in \R^{n\times 1}, g\in SO(n)\right\}$$ and its Lie algebra is
$$\fr(n)=\left\{\bpm 0& 0\\ y & \xi\epm\,\bigg|\, y\in \R^{n\times 1}, \xi\in so(n)\right\}.$$

Given a curve $\g$ in $\R^n$ parametrized by its arc-length, there exists an orthonormal frame $(e_0, \ldots, e_{n-1})$ satisfying $e_0=\g_x$ and $(e_1, \ldots, e_{n-1})$ is a parallel  orthonormal frame for $\nu(\g)$, i.e., 
$$(e_0, \ldots, e_{n-1})_x=(e_0, \ldots, e_{n-1})\bpm 0& -k_1 &-k_2 &\cdots &-k_{n-1}\\ k_1&0&\cdot&&0\\ \cdot & &&&\\k_{n-1}&0&& &0\epm.$$
The orthonormal frame $(e_0, \ldots, e_{n-1})$ is called a {\it p-frame\/} and $k_i$ the {\it principal curvature along $e_i$\/} for $\g$ for $1\leq i\leq n-1$.  Since two p-frames are differed by a constant in $SO(n-1)$, $(\g, e_0, \ldots, e_{n-1})$ is a  $(SO(n-1), Y_n)$-moving frame for $\g\in \calM$, where 
$$\calM=\{\g:\R\to \R^n\n ||\g_x||=1\},$$ and $Y_n$ is the following affine subspace of $\fr(n)$,   
$$ Y_n=e_{21} + \oplus_{i=3}^{n+1} \R(e_{i2}- e_{2i}).$$ 
 \eeg
 
 \beg[{\bf Frenet frame for curves in $\R^3$}]\

The classical Frenet frame $(e_0, \rn, \rb)$ for $\g$ in
$$\calM_s= \{\g:\R\to \R^3\n ||\g_x||=1, \, ||\g_{xx}||>0\}$$
satisfies the following Frenet equation, 
$$(e_0, \rn, \rb)_x=(e_0, \rn, \rb)\bpm 0& -k &0\\ k& 0& -\tau\\ 0& \tau &0\epm,$$
wheat $k, \tau$ are the {\it curvature\/} and {\it torsion\/} of $\g$.  So $(\g, e_0, \rn, \rb)$ is a $(e, V)$-moving frame for $\g\in \calM_s$, where $V= e_{21} + \R(e_{32}-e_{23}) + \R (e_{43}-e_{34})$.  
\eeg

 \beg[{\bf periodic h-frame\/}]\label{es}\

\ni
For $c\in \R$, let $R(c)=\bpm \cos c& -\sin c\\ \sin c& \cos c\epm$ denote the rotation of $\R^2$ by angle $c$. Let $c_0\in \R$ be a constant, and
$$\calM_{c_0}=\{\g:S^1\to \R^3\n ||\g_x||=1, \, {\rm the\, normal\, holonomy\, of \,} \g \, {\rm is\,} R(-2\pi c_0)\}.$$
If  $(e_0, e_1, e_2)$ is a p-frame  along $\g\in \calM_{c_0}$, then
  $$(e_0, e_1, e_2)(2\pi)= (e_0, e_1, e_2)(0)\diag(1, R(-2\pi c_0)).$$ 
 Let $(v_1(x), v_2(x))$ be the orthonormal normal frame obtained by rotating $(e_1(x), e_2(x))$ by $c_0x$. Then the new frame 
 $$\ti g(x)= (e_0, v_1, v_2)(x)= (e_0, e_1, e_2)(x) \bpm 1&0\\ 0& R(c_0x)\epm$$
  is periodic in $x$. Moreover, 
$$\ti g^{-1}\ti g_x= \bpm  0 & -\ti k_1& -\ti k_2\\ \ti k_1 &0&- c_0\\ \ti k_2 & c_0 &0\epm$$ 
and $(\ti k_1+ i\ti k_2)(x,t)= e^{-c_0x} (k_1+ ik_2)(x,t)$ are periodic. We call $\ti g=(e_0, v_1, v_2)$ a {\it periodic $h$-frame\/} along $\g$. This is a $(SO(2), V)$-moving frame for $\calM_{c_0}$, where 
$$V= e_{21} +c_0 e_{43} +\R (e_{32}-e_{23}) +\R (e_{42}- e_{24}).
$$
\eeg

\beg[{\bf Central affine moving frame}]\label{ft} (cf. \cite{UP95}, \cite{CIM13})\

The group $SL(n,\R)$ acts on $\rc$ transitively by $g\cdot y= gy$.  
If $\g:\R\to \rc$ satisfies 
$\det(\g, \g_s, \ldots, \g_s^{(n-1)}) >0$ for all $ s\in \R$,
then we can change to parameter $x=x(s)$ such that $\det(\g, \g_x, \ldots, \g_x^{(n-1)})=1$. Such parameter $x$ is called the {\it central affine arc-length parameter\/}.  Let 
$$\calM_n(\R)= \{ \g\in C^\infty(\R, \rc)\n \det(\g, \ldots, \g_x^{(n-1)})=1\}.$$
Give $\g\in \calM_n(\R)$, take the $x$-derivative of $\det(\g, \ldots, \g_x^{(n-1)})=1$ to get 
$$\det(\g, \g_x, \ldots, \g_x^{(n-2)}, \g_x^{(n)})=0.$$
This implies that there exist $u_1, \ldots, u_{n-1}$ such that 
 $$\g_x^{(n)}= u_1 \g+ \ldots + u_{n-1} \g_x^{(n-2)}.$$    
Hence we have
$$g_x=g\bpm 0&0&0&\cdot &\cdot & u_1\\ 1&0&0&\cdot &\cdot & u_2\\ 0&1&0 &\cdot &\cdot &\cdot\\ 0&&&&& \cdot\\ 0&&&1&0&u_{n-1}\\ 0&&&&1&0\epm,$$
where $g=(\g, \g_x, \ldots, \g_x^{(n-1)})$.  
The map $g$ is called the {\it central affine moving frame\/} along $\g$ and $u_i$ the {\it $i$-th central affine curvature\/} of $\g$ for $1\leq i\leq n-1$.  The central affine moving frame $g$ along $\g\in \calM$ is the $(e, V)$-moving frame along $\g$, where $V= b+ \oplus_{i=1}^{n-1} \R e_{in}$ and $b=\sum_{i=1}^{n-1} e_{i+1, i}$. \eeg

\beg[{\bf Adjoint moving frame}]\label{fs}\

Let $G$ be a semi-simple Lie group, $\calG$ the Lie algebra of $G$ equipped with a non-degenerate ad-invariant bi-linear form $\li\, , \ri$, $a\in \calG$, and 
$$M=\{gag^{-1}\n g\in G\}$$
the Adjoint $G$-orbit at $a$ in $\calG$. Let $G_a$ denote the stabilizer of $a$, $\calG_a^\perp$ the orthogonal complement of the Lie algebra $\calG_a$ of $G_a$ with respect to the bi-linear form $\li\, ,\ri$. It was proved in \cite{TerTho01} that given $\g\in C^\infty(\R, M)$, there is a $g:\R\to G$ such that $\g(x)= g(x)a g(x)^{-1}$ and $g^{-1}g_x\in \calG_a^\perp$. Moreover, if $g_1:\R\to G$ is another such map then there exists a constant $c\in G_a$ such that $g_1= gc$.  Hence $g$ is a $(G_a, \calG_a^\perp)$-moving frame along $\g$ in $C^\infty(\R, M)$ and the $\calG_a^\perp$-valued map $g^{-1}g_x$ is the differential invariants for $\g$ defined by $g$. We call such $g$ an {\it Adjoint moving frame\/}.\eeg

\brem\

\ben
\item[(a)] The group \Diff$(\R)$ of diffeomorphisms of $\R$ acts on $C^\infty(\R, M)$ by $f\cdot\g=\g\circ f$.  This action leaves the space $\calI(\R, M)$ of immersions in $C^\infty(\R, M)$ invariant. Two immersed curves lies in the same \Diff$(\R)$-orbit have the same geometry but with different parametrization. 
\item[(b)]  If \Diff$(\R)\cdot \calM$ is not an open subset of $C^\infty(\R, M)$, then a curve flow on $\calM$ is a curve flow on $M$ with constraints. 
\item[(c)] Suppose $\Diff(R)\cdot \calM$ is open in  $C^\infty(\R, M)$. If  the $(H,V)$-moving frames are for curves with a special parameter then $\dim(V)$ is equal to $\dim(M)-1$, otherwise $\dim(V)=\dim(M)$. 
\een
\erem

The notion of the $(H,V)$-moving frame is closely related to the concept of {\it $H$-slice on a $G$-space\/} $X$ defined by Palais. 

\bdefn(\cite{Pal60})
Let $X$ be a $G$-space, and $H$ a closed subgroup of $G$.  A submanifold $S$ of $X$ is a {\it $H$-slice\/} if $S$ satisfies
\ben
\item $G\cdot S$ is open in $X$,
\item $H\cdot S\subset S$,
\item if $(g\cdot S)\cap S\not=\emptyset$ then $g\in H$.
\een
\edefn 

\bthm Let $M=G\cdot p_0= G/K$, where $K= G_{p_0}$.  Let $G$ act on $C^\infty(\R, M)$ by $(c\cdot\g)(x)= c\cdot\g(x)$, and $C^\infty(\R, K)$ act on $C^\infty(\R, \calG)$ by the gauge transformations, 
$$k\ast u= k^{-1}uk + k^{-1}k_x.$$
Then $\Phi([\g])= [g^{-1}g_x]$ defines a bijection 
$$\Phi: C^\infty(\R, M)/G\to C^\infty(\R, \calG)/C^\infty(\R, K),$$
where $[\g]$ is the $G$-orbit at $\g$ and $[u]$ is the $C^\infty(\R,K)$-gauge orbit of $u$. Moreover, a subset $\calM$ of $C^\infty(\R,M)$ admits  $(H,V)$-moving frames if and only if $C^\infty(\R, V)$ is a $H$-slice of the gauge action of $C^\infty(\R, K)$ on $\ti \calM=\{g^{-1}g_x\in C^\infty(\R, G)\n g\cdot p_0\in \calM\}$.
\ethm

\begin{proof}
Given $\g\in C^\infty(\R, M)$, there exists $g:\R\to G$ such that $\g= g\cdot p_0$. Moreover,  $g_1\in C^\infty(\R, G)$ such that $\g= g_1\cdot p_0$ if and only if there exists $k\in C^\infty(\R, G_{p_0})$ such that $g_1= gk$.  So $\g\mapsto [g^{-1}g_x]$ is a well-defined map from $C^\infty(\R, M)$ to $C^\infty(\R, \calG)/C^\infty(\R, K)$.  Note that if $c\in G$ is a constant, then $(cg)^{-1}(cg)_x= g^{-1}g_x$. This proves that $\Phi([\g])= [g^{-1}g_x]$ is well-defined.  Given $u:\R\to \calG$, we can solve $g\in C^\infty(\R, G)$ such that $u=g^{-1}g_x$.  So $\Phi(g\cdot p_0)=u$.  This shows that $\Phi$ is a bijection. 
 The second part of the theorem follows from the definitions.  
\end{proof}

 Fels and Olver gave a general theory and a systematic method of constructing group based moving frames in \cite{FO99}. Mar\'i Beffe used these group based moving frames and bi-Hamiltonian structures to study various integrable curve flows in \cite{MB08a}, \cite{MB08b}, \cite{MB09a}, \cite{MB10a}.  We will see in later sections that our notion of $(H,V)$-moving frames 
makes the relation between integrable geometric curve flows and soliton equations more transparent.

\bs

\section{The $2\times 2$ AKNS hierarchy}\label{bz}

In this section, we review some known properties of the $2\times 2$ AKNS hierarchy and its various restrictions (cf. \cite{AblCla91}, \cite{TerUhl98}).

\ss
\bsub{\bf The $2\times 2$ AKNS hierarchy or the $SL(2,\C)$-hierarchy}\

Let 
$$a=\diag(i, -i).$$ 
It can be checked that given $u=\bpm 0& q\\ r& 0\epm:\R\to sl(2,\C)$ there is a unique 
$$Q(u,\l)= a\l+ Q_0(u) + Q_{-1}(u)\l^{-1} + \cdot$$
satisfying $Q_0(u)=u$,
\beq\label{bk}
\bca
[\p_x+ a\l + u, Q(u,\l)]=0,\\
Q(u,\l)^2= -\l^2\I_2,
\eca
\eeq
where $\I_2$ is the $2\times 2$ identity matrix.
Moreover, entries of $Q_{-j}(u)$ are differential polynomials in $q$ and $r$ (i.e., polynomials in $q, r$ and their $x$-derivatives. It follows from \eqref{bk} that we have the recursive formula
\beq\label{df}
(Q_{-j}(u))_x + [u, Q_{-j}(u)]= [Q_{-(j+1)}(u), a].
\eeq 
In fact,  the $Q_j(u)$'s can be computed directly from \eqref{bk}.  For example, 
\beq\label{bl}
Q_{-1}(u)=\frac{i}{2}\bpm qr & q_x\\ -r_x & -qr\epm, \quad Q_{-2}(u)=-\frac{1}{4}\bpm q_xr - qr_x& q_{xx}- 2q^2r\\ r_{xx}- 2qr^2& -q_xr+ qr_x\epm.
\eeq

Let  $$V=\C e_{12} + \C e_{21}.$$
The {\it $j$-th flow in the $2\times 2$ AKNS hierarchy\/} is the following flow on $C^\infty(\R, V)$,
\beq\label{bm}
u_t=[\p_x+ u, Q_{-(j-1)}(u)]= [Q_{-j}(u), a].
\eeq  
For example,  the first three flows in the $SL(2,\R)$-hierarchy are 
 \begin{align*}
 &q_t= q_x, \quad r_t= r_x,\\
 & q_t= \frac{i}{2} (q_{xx} - 2q^2 r),\quad r_t= -\frac{i}{2} (r_{xx} - 2 qr^2),\\
 & q_t= -\frac{1}{4} (q_{xxx} - 6qrq_x),\quad r_t= -\frac{1}{4}(r_{xxx} - 6qr r_x).
 \end{align*}

It follows from \eqref{bk} that $u:\R^2\to V$ is a solution of \eqref{bm} if and only if
\beq\label{dc}
\o_j=(a\l + u) \rd x+\left( \sum_{-(j-1)\leq i\leq 1} Q_i(u)\l^{j-1+i}\right) \rd t
\eeq
is a flat $sl(2,\C)$-valued connection one form on the $(x,t)$-plane for all complex parameter $\l$, where $Q_1(u)= a$ and $Q_0(u)=u$.  The connection $1$-form $\o_j$ is the {\it Lax pair of the solution $u$ of the $j$-th flow\/} \eqref{bm}. 
\esub

\ss
\bsub {\bf The $SU(2)$-hierarchy}\

Let $V_1$ denote the following linear subspace of $V=\C e_{12} +\C e_{21}$:
\beq\label{en}
V_1=\{ qe_{12}-\bar q e_{21}\n q\in \C\}.
\eeq
 The $j$-th flow \eqref{bm} leaves $C^\infty(\R, V_1)$-invariant, i.e., leaves $r=-\bar q$ invariant. Moreover, $Q_{-j}(u)\in su(2)$ for $u=\bpm 0& q\\ -\bar q&0\epm$.
 The restrictions of the flows \eqref{bm} to $C^\infty(\R, V_1)$ is called the {\it $SU(2)$-hierarchy\/}.
For example, for $u=\bpm 0& q\\ -\bar q &0\epm$, we have
\begin{align}
&Q_{-1}(u)=\frac{i}{2}\bpm -|q|^2& q_x\\ \bar q_x & |q|^2\epm, \label{db1}\\
& Q_{-2}(u)= \frac{1}{4}\bpm \bar q q_x- q\bar q_x& -(q_{xx} + 2|q|^2 q)\\ \bar q_{xx} + 2|q|^2 \bar q_x & -\bar q q_x+ q\bar q_x\epm.\label{db2}
\end{align} 
So the first three flows in the $SU(2)$-hierarchy are
\begin{align}
&q_t= q_x,\\
& q_t= \frac{i}{2}(q_{xx}+ 2|q|^2 q), \label{bo2}\\
& q_t= -\frac{1}{4}(q_{xxx} +6|q|^2 q_x). \label{bo3}
\end{align}
Note that the second flow is the focusing NLS and the third equation is the {\it complex modified KdV\/}. 
 The Lax pair $\o_j$ (defined by \eqref{dc}) for solution $u=\bpm 0& q\\ -\bar q &0\epm$ of the $j$-th flow in the $SU(2)$-hierarchy  satisfies the {\it $su(2)$-reality condition\/}
\beq\label{cb}
\o_j(x,t,\bar\l)^*+ \o_j(x,t,\l)=0.
\eeq
We call a solution $E$ of 
$$E^{-1}\rd E:= E^{-1}(E_x \rd x+ E_t \rd t)= \o_j$$ a {\it frame\/} of the $j$-th flow in $SU(2)$-hierarchy if $E$ satisfies the {\it $SU(2)$-reality condition\/},
\beq\label{cd}
E(x,t,\bar\l)^*E(x,t,\l)=\I_2.
\eeq
\esub

\ss
\bsub {\bf The $SU(1,1)$-hierarchy}\

Let 
\beq\label{en1}
V_2=\{ qe_{12}+ \bar q e_{21}\n q\in \C\}.
\eeq
The $j$-th flow \eqref{bm} leaves $C^\infty(\R, V_2)$ invariant, i.e., leaves $r= \bar q$ invariant. Moreover, $Q_{-j}(u)\in su(1,1)$ for $u=\bpm 0& q\\ \bar q&0\epm$. The restriction of the $2\times 2$ AKNS hierarchy to $C^\infty(\R, V_2)$ is called the $SU(1,1)$-hierarchy. 
For $u=\bpm 0& q\\ \bar q &0\epm$, we have 
\begin{align}
& Q_{-1}(u)=\frac{i}{2}\bpm |q|^2& q_x\\ -\bar q_x & -|q|^2\epm, \label{aq} \\
& Q_{-2}(u)= \frac{1}{4}\bpm -\bar q q_x+q\bar q_x& -(q_{xx} - 2|q|^2 q)\\
-( \bar q_{xx} - 2|q|^2 \bar q)& \bar q q_x-q\bar q_x\epm. \label{aqa}
\end{align}
The first three flows in the $SU(1,1)$-hierarchy are
\begin{align}
&q_t= q_x,\\
& q_t=\frac{i}{2}(q_{xx} - 2|q|^2 q), \label{bp2}\\
&q_t=  -\frac{1}{4} (q_{xxx}- 6|q|^2 q_x). \label{bp3}
\end{align}
Note that the second flow is the {\it defocusing NLS\/}.  
The Lax pair $\o_j$ for a solution $u=\bpm 0& q\\ \bar q&0\epm$ of the $j$-th flow of the $SU(1,1)$-hierarchy is given by the same formula \eqref{dc} and  $\o_j$ satisfies the $su(1,1)$-reality condition,
\beq\label{ce}
\o_j(x,t,\bar\l)^*J_2+ J_2\o(x,t,\l)=0, \quad J=\diag(1,-1).
\eeq
We call a solution $E$ of $E^{-1}\rd E= \o_j$ a {\it frame\/} if $E$ satisfies the {\it $SU(1,1)$-reality condition\/},
\beq\label{cf}
E(x,t,\bar\l)^*JE(x,t,\l)=J.
\eeq
\esub

\ss
\bsub {\bf The $SL(2,\R)$-hierarchy}
 
 Let $$a= \diag(1,-1), \quad u=\bpm 0& q\\ r& 0\epm,$$
 and $Q(u,\l)= a\l + u+ Q_{-1}(u)\l^{-1} +\cdots$
  the solution of 
 \beq\label{bka}
 \bca [\p_x+ a\l + u, Q(u,\l)]=0,\\ Q(u,\l)^2=\l^2\I_2.\eca
 \eeq
  Then
 \begin{align}
 Q_{-1}(u)&=\frac{1}{2}  \bpm -qr& -q_x\\ r_x & qr\epm,\label{bv1} \\
 Q_{-2}(u)&= \frac{1}{4}\bpm q_x r- qr_x & q_{xx}- 2q^2r\\ r_{xx} - 2qr^2 & -q_xr + qr_x\epm.\label{bv2}
 \end{align}
  Let 
 \beq\label{fk}
 V_\R= \R e_{12}+ \R e_{21}.
 \eeq
 The {\it $j$-th flow on $C^\infty(\R, V_\R)$ in the $SL(2,\R)$-hierarchy\/} is given by the same formula \eqref{bm} with these new $a$ and $Q_j$'s.
For example,  the first three flows in the $SL(2,\R)$-hierarchy for $u=\bpm 0& q\\ r& 0\epm$ are 
 \begin{align*}
 &q_t= q_x,\\
 & q_t= -\frac{1}{2} (q_{xx} - 2q^2 r),\quad r_t= \frac{1}{2} (r_{xx} - 2 qr^2), \\
 & q_t= \frac{1}{4} (q_{xxx} - 6qrq_x),\quad r_t= \frac{1}{4}(r_{xxx} - 6qr r_x).
 \end{align*}
  The Lax pair $\o_j$ of the $j$-th flow in the $SL(2,\R)$-hierarchy, which is given by the same formula \ref{dc}, satisfies the $sl(2,\R)$-reality condition,
 \beq\label{cg}
 \overline{\o_j(x,t,\bar\l)}=\o_j(x,t,\l).
 \eeq
A solution of $E^{-1}\rd E= \o_j$ is a {\it frame\/} of the solution $u$ of the $j$-th flow in the $SL(2,\R)$-hierarchy if $E$ satisfies the {\it $SL(2,\R)$-reality condition\/},
 \beq\label{ch}
\overline{ E(x,t,\bar\l)}=E(x,t,\l).
 \eeq
\esub

\bsub {\bf The KdV hierarchy}\

The $(2j+1)$-th flow in the $SL(2,\R)$-hierarchy leaves $r=1$ invariant and the third flow becomes the KdV, $q_t= \frac{1}{4}(q_{xx} - 6 qq_x)$. 
The Lax pair of the KdV is
\beq\label{fq}
\o= (a\l + u)\rd x + (a\l^3 +
u\l^2 + Q_{-1}(u)\l + Q_{-2}(u))\rd t,
\eeq
where 
\begin{align*}
& a= \diag(1,-1), \quad u=\bpm 0& q\\ 1& 0\epm, \\
& Q_{-1}(u)= \bpm -\frac{q}{2}& -\frac{q_x}{
2}\\ 0 & \frac{q}{ 2}\epm,\quad Q_{-2}(u)= \bpm \frac{q_x}{4}& \frac{q_{xx} -
2q^2}{ 4}\\ -\frac{q}{ 2}& -\frac{q_x}{ 4}\epm.
\end{align*}
The reality condition for the KdV hierarchy is a little bit more complicated. It was noted in \cite{TU00} that the Lax pair \eqref{fq} satisfies the {\it KdV reality condition\/},
\beq\label{ex}
\bca \overline{A(\bar \l)}=A(\l), \\
\phi(\l)^{-1}A(\l)\phi(\l) =
\phi(-\l)^{-1} A(-\l) \phi(-\l),
\eca\eeq
where $$\phi(\l)=\bpm 1& \l \\0& 1\epm.$$
A frame of a solution $q$ of the KdV is the solution $E(x,t,\l)$ of the following system,
$$\bca E_x= E(a\l+u), \quad E_t= E(a\l^3+ u\l^2+ Q_{-1}(u)\l + Q_{-2}(u)),\\
\overline{E(x,t,\bar\l)}=E(\l), \quad \phi^{-1}(\l) E(x,t,\l) \phi(\l) \,\, {\rm is\, even\, in \, \,} \l.\eca
$$
\esub

\brem
The reality condition of the Lax pair of a soliton equation plays an important role in the symmetries of the equation.  It also plays an essential role in the construction of B\"acklund transformations. 
\erem

\ss
\bsub {\bf Bi-Hamiltonian structure} (\cite{DS84}, \cite{Ter97})
 
Let $\calG= su(2), su(1,1)$ and $sl(2,\R)$ for $W= V_1, V_2$ and $V_\R$ defined by \eqref{en}, \eqref{en1} and \eqref{fk} respectively. The {\it gradient\/} of $H:C^\infty(S^1, W)\to \R$ is the unique $\K H(u)\in C^\infty(S^1, W)$ defined by
 $$\rd H_u(v)= \oint \li \K F(u), v\ri \rd x,$$
 where $\li X, Y\ri= -\frac{1}{2} \tr(XY)$. These hierarchies are Hamiltonian with respect to two Poisson structures. The first one is 
 \beq\label{bn1}
 \{F, H\}_1(u)=\oint \li [\K F(u), a], \K H(u)\ri \rd x
 \eeq
To define the second Poisson structure we need the operator  
$$P_u:C^\infty(S^1, W)\to C^\infty(S^1, \calG)$$ defined as follows: Given $u,v\in C^\infty(\R,W)$ there is a unique $\ti v\in C^\infty(\R, \calG)$ such that $\ti v= u+ A a$ and $[\p_x+u, \ti v]\in W$.  Then $P_u(v)=\ti v$.  
For example,  for $u=\bpm 0& q\\ r&0\epm$ and $v=\bpm 0& \xi\\ \eta&0\epm$, we have $P_u(v) = v+ A a$, where $A$ satisfies
$$A_x= -q\eta+ r \xi.$$
The second Poisson structure is 
\beq\label{bn2}
 \{F, H\}_2(u)= \oint \li [\p_x+ u, P_u(\K F(u))], \K H(u)\ri \rd x.
 \eeq 
  The Hamiltonian equation for $H$ with respect to $\{\, ,\}_i$ is 
 $$u_t= (J_i)_u(\K H(u)),$$
 where 
 \begin{align}
 &(J_1)_u(v)= [v, a],\label{bna}\\
 &(J_2)_u(v)= [\p_x+u, P_u(v)].  \label{bnb}
 \end{align}
 Note that $\{\, ,\}_2$ are defined for $H$ satisfying 
$\oint \li u, [v,a]\ri \rd x=0$.  

Let $H_j:C^\infty(S^1, W)\to \R$ be the functional defined by
\beq\label{bq}
H_j(u)=-\frac{1}{j}\oint \li Q_{-j}(u), a\ri \rd x.
\eeq
Then the following are known:
\ben
\item[(i)] The $j$-th flow equation \eqref{bm} is the Hamiltonian equation for $H_j$ with respect to $J_2$ and is the Hamiltonian equation for $H_{j+1}$ with respect to $J_1$. 
\item[(ii)]  Write 
$$Q_j(u)= Y_j(u) + A_j(u) a$$
with $Y_j(u)\in W$ and $A_j(u)\in \R$. Then we have
\begin{align*}
&P_u(Y_{-j}(u))= [Q_{-(j+1)}(u),a],\\
&Y_{-1}(u)= (J_1)_u^{-1}(u_x),\\
&\K H_j(u)= Y_{-(j-1)}(u).
\end{align*}
This gives the well-known fact that the $j$-th flow can be written recursively by the bi-Hamiltonian structures,
$$u_{t_j}= (J_2J_1^{-1})_u^j (u_x).$$
\item[(iii)] $H_j$'s are commuting Hamiltonians with respect to both Poisson structures.  Hence all flows commute and $H_j$'s are conserved quantities for these flows. 
\een 
\esub

\ss
\bsub {\bf The $\R$-action on solutions\/} 

These hierarchies also admit an $\R$-action on solutions. Let $\R$ act on $C^\infty(\R, W)$ by 
$$c\ast u= e^{ca} ue^{-ca},$$
where $a=\diag(i, -i)$ for $W=V_1$ and $V_2$ and $a= \diag(1,-1)$ for $W=V_\R$. 
For $u\in V_1$ or $V_2$, the induced $\R$-action on $q\in C^\infty(\R, V_i)$ is
$$c\ast q= e^{2ic} q.$$
So this gives an $S^1$-action on $C^\infty(\R, V_i)$. 
For $u=qe_{12}+ r e_{21}\in V_\R$, the $\R$-action on $(q,r)$ is
 $$c\ast (q, r)= (e^{2c} q, e^{-2c} r).$$
 
 It follows from \eqref{bk} that 
$$Q(e^{ca} u e^{-ca})= e^{ca}Q(u) e^{-ca}.$$
So we have 
\beq\label{dj}
Q_{-j}(e^{ca} u e^{-ca})= e^{ca}Q_{-j}(u) e^{-ca}.
\eeq
This implies that
\ben
\item if $u$ is a solution of the $j$-th flow then so is $e^{ca} u e^{-ca}$ for any $c\in \R$.  In other words, $S^1$ acts on the space of solutions of the $j$-th flow in the $SU(2)$ and $SU(1,1)$ hierarchies and the group $\R^+=\{\diag(r, r^{-1})\n r\in \R, r\not=0\}$ acts on the space of solutions of the $j$-th flow in the $SL(2,\R)$-hierarchy,
\item the functional $H_j$ defined by \eqref{bq} is invariant under the $\R$-action.  
 \een
 
 Note that if $F:C^\infty(S^1,W)\to \R$ is invariant under the $\R$-action, then $\K F$ is $\R$-equivariant, i.e., 
  $$\K F(e^{ca} u e^{-ca})= e^{ca} \K F(u) e^{-ca}.$$
 Hence we get the following.
 
 \bprop\label{dd}
 If $F, H:C^\infty(S^1, W)\to \R$ are invariant under the $\R$-action $r\ast u= e^{ra}u e^{-ra}$, then $\{F, H\}_i$ is also $\R$-invariant for $i=1,2$. 
  \eprop
   \esub
   
   Note that the Hamiltonian theory discussed above works for the space $\calS(\R, W)$ of  rapidly decaying smooth maps from $\R$ to $W$. 
   
 \bs
\section{The $\ast$-MCF on $\R^3$}\label{ca}

In this section, we 
\ben
\item[(a)] show that the $\ast$-MCF on $\R^3$ preserves arc-length parameter and is the VFE,
\item[(b)] explain the Hasimoto transform for the VFE in terms of p-frames,
\item[(c)] solve the Cauchy problem for the VFE with initial data having rapidly decaying principal curvatures.
\een

Recall that the Hodge star operator on an oriented two dimension inner product space is the rotation of $\frac{\pi}{2}$.  So if $(e_1, e_2)$ is an oriented orthonormal basis then 
$$\ast(e_1)= e_2, \quad \ast(e_2)= -e_1.$$

\bprop  The $\ast$-MCF on a $3$-dimensional Riemannian manifold $(N^3, \rg)$ preserves the arc-length parameter.
\eprop

\begin{proof} 
Suppose $\g(x,t)$ is a solution of the $\ast$-MCF on $N$. Let $s(\cdot,t)$ denote the arc-length parameter of $\g(\cdot, t)$. Then $s_x= ||\g_x|| =\sqrt{\rg(\g_x, \g_x)}$.  Let $(e_0, e_1, e_2)$ be an orthonormal moving frame along $\g$ such that $e_0$ is the unit tangent, and $k_i$ the principal curvature of $\g$ along $e_i$ for $i=1,2$. Then  
\beq\label{bj}
(e_0, e_1, e_2)_s= (e_0, e_1, e_2)\bpm 0& -k_1 &-k_2\\ k_1&0&-\b\\ k_2 &\b &0\epm
\eeq
for some functions $k_1, k_2$ and $\b$. The $\ast$-MCF on $N$ written in terms of this frame is
$$\g_t= \ast_\g (H(\g(\cdot, t)))= k_1 e_2- k_2 e_1.$$  Write $||\g_x||^2=\rg(\g_x,\g_x)$. Then 
\begin{align*}
\frac{1}{2}(\rg(\g_x, \g_x))_t&=\rg(\K_{\frac{\p}{\p t}}\g_x, \g_x)= \rg(\K_{\g_x} \g_{t}, \g_x) \\
&= \rg(\K_{\g_x}(k_1 e_2- k_2 e_1), \g_x)= ||\g_x||^2\rg(\K_{e_0} (k_1e_2- k_2 e_1), e_0).
\end{align*}
It follows from \eqref{bj} that this is zero.
\end{proof}

So we may assume solutions $\g(x,t)$ of the $\ast$-MCF is parametrized by arc-length. 
Since $H(\g)= \g_{xx}$ if $\g$ is parametrized by its arc-length, we obtain the following.

\bprop The $\ast$-MCF on the Euclidean space $\R^3$ parametrized by  arc-length is the {\it vortex filament equation\/} (VFE), 
$$\g_t=\g_x\times\g_{xx}= k{\bf b},$$ where $k(\cdot, t)$ and ${\bf b}(\cdot, t)$ are the curvature and bi-normal of $\g(\cdot, t)$.
\eprop  

Let $c\in \R\bh 0$ be a constant. Then $\g$ is a solution of $\g_t= c\g_x\times \g_{xx}$ if and only if $\ti \g(x,t)= \g(x, t/c)$ is a solution of the VFE. So they are the same equation up to linear change of time coordinate.  To make the relation between the VFE and NLS looks neater, we will call 
\beq\label{cr}
\g_t= \frac{1}{2}\ast_\g H(\g(\cdot, t))= \frac{1}{2} \g_x\times \g_{xx}
\eeq
 the VFE for the rest of the paper. 

Next we state Hasimoto's result in terms of p-frame. Since the proof of this theorem is used for many integrable curve flows, we include a proof here. 

\bthm\label{bb} If $\g:\R^2\to \R^3$ is a solution of the VFE, $\g_t= \frac{1}{2}\g_x\times \g_{xx}$, parametrized by arc-length, then there exists $g=(e_0, e_1, e_2):\R^2\to SO(3)$ satisfying
\ben
\item[(i)] $g(\cdot, t)$ is a $p$-frame along $\g(\cdot, t)$ for each $t$,
\item[(ii)] $q=\frac{1}{2}(k_1+ik_2)$ is a solution of the NLS, $q_t=\frac{i}{2}(q_{xx}+ 2|q|^2 q)$, where $k_i=(e_0)_x\cdot e_i$ is the principal curvatures along $e_i$ for $i=1,2$. 
\item[(iii)] if both $g$ and $\ti g=(e_0, \ti e_1, \ti e_2)$ satisfies (i) and (ii), then there is a constant $c$ such that $\ti g= g\diag(1, R(c))$ and $\ti q=\frac{1}{2}(\ti k_1+i\ti k_1)= e^{-ic} q$ is a solution of the NLS, where $\ti k_i(\cdot, t)$ is  the principal curvature of $\g(\cdot, t)$ along $\ti e_i$ and $R(c)=\bpm \cos c& -\sin c\\ \sin c & \cos c\epm$.
\een 
\ethm

\begin{proof}
Let $h=(e_0, v_1, v_2):\R^2\to SO(3)$ such that $h(\cdot, t)$ is a p-frame for $\g(\cdot, t)$ for all $t$ and $r_1, r_2$ the principal curvatures along $v_1, v_2$ respectively, i.e., 
$$A:= h^{-1}h_x= \bpm 0& -r_1 & -r_2\\ r_1&0 &0\\ r_2&0&0\epm.$$
Let $$B= (b_{ij}):= h^{-1}h_t.$$
 A direct computation implies that 
$$(e_0)_x= (\g_x)_t= (\g_t)_x= \frac{1}{2}(-r_2 v_1+ r_1 v_2)_x=\frac{1}{2} (-(r_2)_x v_1 + (r_1)_x v_2.$$
So $b_{21}= -\frac{1}{2} (r_2)_x$ and $b_{31}= \frac{1}{2} (r_1)_x$.  
Since $h^{-1}h_x= A$ and $h^{-1}h_t= B$, we have 
$$A_t= B_x+ [A, B].$$
Compare the 32-th entry of the above equation to get
$$(b_{32})_x= -\frac{1}{4} (r_1^2+ r_2^2)_x.$$
Hence there exists $\o:\R\to \R$ such that $b_{32}= -\frac{1}{4} (r_1^2+ r_2^2)+ \o(t)$.  Let $(e_1, e_2)(x,t)= (v_1 v_2)(x,t)R(\rho(t))$, and $g=(e_0, e_1, e_2)$, where $R(\rho(t))$ is the rotation of angle $\rho(t)$ and $\rho'(t)=-\o(t)$. Then we have
\beq\label{di}
\bca g^{-1}g_x= \bpm 0& -k_1& -k_2\\ k_1& 0&0\\ k_2&0&0\epm, \\
 g^{-1}g_t= \bpm 0 & \frac{(k_2)_x}{2} & -\frac{(k_1)_x}{2} \\ -\frac{(k_2)_x}{2} &0 &\frac{k_1^2+k_2^2}{4}\\
\frac{(k_1)_x}{2} &-\frac{k_1^2+k_2^2}{4} & 0\epm,\eca
\eeq
where $k_1+ ik_2= e^{-i\rho(t)} (r_1+ r_2)$.  Let $\ti A= g^{-1}g_x$ and $\ti B= g^{-1}g_t$. Then $\ti A_t= (\ti B)_x+[\ti A, \ti B]$.  Equate the 21-th and 31-th entries of this equation to see that $q=q_1+iq_2= \frac{1}{2} (k_1+ ik_2)$ is a solution of the NLS. 

Statement (iii) follows from the construction of $g$. 
\end{proof}

\bcor Let $\g:\R^2\to \R^3$ be a solution of $\g_t=\frac{1}{2}\g_x\times \g_{xx}$ parametrized by the arc-length, $g=(e_0, e_1, e_2):\R^2\to SO(3)$ such that $g(\cdot, t)$ is a p-frame along $\g(\cdot, t)$, and $r_i(\cdot, t)$ the principal curvature with respect to $e_i(\cdot, t)$ for $i=1,2$. Then there exists a smooth function $\o:\R\to \R$ such that $q(x,t)= \frac{1}{2}e^{i\o(t)} (r_1+ ir_2)$ is a solution of the NLS. 
\ecor

We have seen in section \ref{bz} that  $q$ is a solution of the NLS  if and only if 
$$\o_2(\cdot, \cdot, \l)= (a\l+ u) \rd x + (a\l^2 + u\l + Q_{-1}(u))\rd t,$$
is a flat connection $1$-form on the $(x,t)$-plane for all complex parameter $\l$, where $a=\diag(i,-i)$, $u=qe_{12}-\bar q e_{21}$, and $Q_{-1}(u)$ is defined by \eqref{db1}. Recall that 
a {\it frame\/} of the solution $u$ of the NLS is a solution $E(x,t,\l)$ of $E^{-1}\rd E= \o_2$ satisfying the $SU(2)$-reality condition,
\beq\label{as}
E(\cdots, \bar\l)^*E(x,t, \l)=\I.
\eeq

The following proposition follows from the fact that $\tr(\o_2)=0$ and $E$ satisfies the $SU(2)$-reality condition.

\bprop\label{eq}
If $E(x,t,\l)$ is a frame for the solution $q$ of the NLS, then 
\ben
\item $\det(E(x,t,\l))$ is independent of $x,t$,
\item $E(x,t,r)\in U(2)$ for $r\in \R$.
\een
\eprop
 
 Next we identify the Euclidean $\R^3$ as $su(2)$. Let $su(2)$ be equipped with inner product 
 $$(X,Y)=-\frac{1}{2} \tr(XY).$$ Let 
 \beq\label{cq} 
 a=\diag(i, -i), \quad b= \bpm 0& 1\\ -1 &0\epm, \quad c= \bpm 0& i\\ i&0\epm.
 \eeq
 Then $\d=(a, -c, b)$ is an ordered orthonormal basis for $su(2)$. 
Note that
 $$[a,b]= 2c, \quad [b,c]=2 a, \quad [c,a]=2 b.$$  We identify $su(2)$ as the Euclidean $\R^3$ via
 \beq\label{ep}
 xa+yb+zc=\bpm ix & y+iz\\ -y+iz & -ix\epm \quad \rightarrow \bpm x\\ y\\ z\epm.
 \eeq
  
 Let $\Ad: SU(2)\to SO(su(2))$ denote the Adjoint representation, i.e., 
 $$\Ad(g)(\xi)= g\xi g^{-1}.$$   Let 
 \beq\label{cqa}
 \Ad(g)[\d]=(gag^{-1}, -gcg^{-1}, gbg^{-1}),
 \eeq
 where $a, b, c$ are as in \eqref{cq}. 
 Then $\Ad(g)[\d]$ is in $SO(3)$ (here we identify $su(2)$ as $\R^3$ by \eqref{ep}).    
  
Below is the standard Pohlmeyer-Sym construction in soliton theory (cf. \cite{Poh76}, \cite{Sym85}). Since we will use this often, we include the elementary proof here.
 
 \bthm\label{bda} Let $E(x,t,\l)$ be a frame of a solution $q:\R^2\to \C$ of the NLS such that $E\in SL(2,\C)$. Then $$\a= \frac{\p E}{\p \l}E^{-1}\,\big|_{\l=0}$$ lies in $su(2)$ and is a solution of the VFE, $\g= \frac{1}{2} \g_x\times \g_{xx}$, parametrized by the arc-length. Moreover, let $\phi(x,t)= E(x,t,0)$. Then $g=\Ad(\phi)[\d]$ defined by \eqref{cqa} satisfies \eqref{di} with $k_1+ik_2= 2q$.  
  \ethm
  
  \begin{proof} Since $E$ satisfies the $SU(2)$-reality condition, $E(x,t,\l)\in U(2)\cap SL(2,\C)= SU(2)$ for all $\l\in \R$.  So $\a\in su(2)$.   Let $\phi(x,t)= E(x,t,0)$. Then $\phi\in SU(2)$ and
  \begin{align}
 & \phi^{-1}\phi_x= u=\bpm 0& q\\ -\bar q &0\epm,\label{dy1},\\ & \phi^{-1}\phi_t= Q_{-1}(u)= \frac{i}{2}\bpm -|q|^2 & q_x\\ \bar q_x & |q|^2\epm,\label{dy2}
 \end{align}
 where $u=qe_{12}- \bar q e_{21}$.  Use $E^{-1}E_x= a\l +u$, $E^{-1}E_t= a\l^2 + u\l + Q_{-1}(u)$ and a direct computation to see that
 \beq\label{dz}
 \bca \a_x= \phi a\phi^{-1}, \\ 
\a_t = \phi u\phi^{-1}
\eca 
\eeq
  where $q= q_1+ iq_2$.  
  
  So $\a(\cdot, t)$ is parametrized by the arc-length.  Let 
  $$e_0=\phi a\phi^{-1}, \quad e_1=-\phi c \phi^{-1}, \quad e_2= \phi b \phi^{-1}.$$  Since $\phi\in SU(2)$ and $(a, -c, b)$ is an orthonormal basis of $su(2)$, $(e_0, e_1, e_2)\in SO(3)$.  Write $q=q_1+ iq_2$. Compute directly to get
  \begin{align*}
  (e_0)_x&= \phi [\phi^{-1}\phi_x, a] \phi^{-1}= \phi [u, a]\phi^{-1}\\ &= 2 q_2 \phi b\phi^{-1}- 2 q_1\phi c\phi^{-1}= 2q_1 e_1+ 2 q_2 e_2,\\
  (e_1)_x&= -\phi [u,c]\phi^{-1}= -2q_1 e_0.
  \end{align*}
  This shows that 
  $$(e_0, e_1, e_2)_x= \bpm 0& - 2q_1& -2q_2\\ 2q_1&0&0\\ 2q_2 &0&0\epm,$$
  i.e., $(e_0, e_1, e_2)(\cdot, t)$ is a p-frame along $\g(\cdot, t)$ and $2q_1, 2q_2$ are the principal curvatures along the parallel normal $e_1, e_2$ respectively. 
  Use \eqref{dy2} and a similar computation to get
    $$(e_0, e_1, e_2)_t= (e_0, e_1, e_2)\bpm 0 & (q_2)_x & -(q_1)_x\\ -(q_2)_x & 0 & |q|^2\\ (q_1)_x & -|q|^2 & 0\epm.$$ 
  
  The second equation of \eqref{dz} implies that $\a_t = -q_2 e_1 + q_1 e_2= \frac{1}{2}\ast_\g H(\g(\cdot, t))$. 
 \end{proof}
 
 \brem The condition $E\in SL(2,\C)$ in Theorem \ref{bad} is not essential. To see this, let $E$ be a frame of the solution $q$ of the NLS, and $\a=\frac{\p E}{\p\l}E^{-1}|_{\l=0}$.  By Proposition \ref{eq}, $\det(E(x,t,\l))$ is independent of $x, t$. So $\tr(\a)$ is a constant.  Hence $\g:=\a-\frac{1}{2}\tr(\a)\I_2\in su(2)$ and is a solution of $\g_t= \frac{1}{2}\g_x\times \g_{xx}$.  
  \erem
  
  If $E_1, E_2$ are two frames for the solution $q$ of the NLS satisfying $E_i(0,0, \l)\in SU(2)$ for $\l\in \R$, then there is $f(\l)\in SL(2,\C)$ satisfies the $SU(2)$-reality condition such that $E_2(x,t,\l)= f(\l) E_1(x,t,\l)$ and $f(\l)\in SU(2)$ for $\l\in \R$.  So we obtain the following.
       
   \bprop\label{er}
   Let $E_1, E_2\in SL(2,\C)$ be frames of the solution $q$ of the NLS, and $\g_i=\frac{\p E_i}{\p \l} E_i^{-1}|_{\l=0}$ the solution of $\g_t= \frac{1}{2} \g_x\times \g_{xx}$ constructed from $E_i$ for $i=1,2$. Then there is a rigid motion $\phi$ of $\R^3$ such that $\g_2= \phi \circ \g_1$.
   \eprop
  
Given a solution $\g$ of the VFE  \eqref{cr}, let $q$ be a solution associated to $\g$ in Theorem \ref{bb}.  Use Theorem \ref{bda} to construct solutions $\ti \g$ of \eqref{cr} from $q$.  The following theorem gives the relation between $\ti \g$ and $\g$.

\bthm\label{bd}
Let $\g:\R^2\to \R^3$ be a solution of the VFE, $\g_t= \frac{1}{2} \g_x\times \g_{xx}$ parametrized by the arc-length, $g:\R^2\to SO(3)$, $k_1, k_2$ principal curvatures, and $q=\frac{1}{2}(k_1+i k_2)$ the solution of the NLS  as in Theorem \ref{bb}. Choose $\phi_0\in SU(2)$ such that $g(0,0)=\Ad(\phi_0)[\d]$. Let $E$ be the frame of $q$ satisfying $E(0,0,\l)=\phi_0$, and $\a=\frac{\p E}{\p\l}E^{-1}\n_{\l=0}$. Then $\g= \a+\g(0,0)$ and $g= \Ad(\phi)[\d]$, where $\phi(x,t)=E(x,t,0)$.  
\ethm

\begin{proof} Write $g=(e_0, e_1, e_2)$, $\phi(x,t)= E(x,t,0)$, and 
$$h=(\ti e_0, \ti e_1,\ti e_2)=\Ad(\phi)[\d]=(\phi a \phi^{-1}, -\phi c\phi^{-1}, \phi b\phi^{-1}).$$ By Theorem \ref{bda}, $g$ and $h$ satisfies the same linear system \eqref{di} and have the same initial data. So $g=h$, which implies that $\ti e_i= e_i$ for $0\leq i\leq 2$.  But $\g_x= e_0= \ti e_0=\a_x$.  We also have 
$\g_t= \frac{1}{2} (k_1 e_2- k_2 e_1)=\a_t$. Since $E(0,0,\l)=\phi_0$ is independent of $\l$,  $\a(0,0)=0$. Therefore $\g= \a+ \g(0,0)$. 
\end{proof} 

It is clear that the VFE $\g_t=\frac{1}{2}\g_x\times \g_{xx}$ is invariant under the group $\fR(3)$ of rigid motions of $\R^3$. We have seen that the space of solutions of the NLS is invariant under the $S^1$-action.  As a consequence of Theorems \ref{bda}, \ref{bd}, and Proposition \ref{er}, we have:

 \bthm\label{cx} Let $\calC=\{\g\in C^\infty(\R^2, \R^3)\n ||\g_x||=1, \g_x= \frac{1}{2}\g_x\times \g_{xx}\}$, $\fR(3)$ the rigid motion group of $\R^3$, and $\calN$ the space of smooth solutions of the NLS, and $\Psi:\calC/\fR(3)\to \calN/S^1$ defined by $\Psi([\g])= [q]$, where $[\g]$ is the $\fR(3)$-orbit of $\g$ and  $[q]$ is the $S^1$-orbit of  a solution $q$ constructed from $\g$ in Theorem \ref{bb}. Then
$\Psi$ is well-defined and is a bijection.
 \ethm

Theorem \ref{bd} also implies the following.

\bthm[Cauchy problem for VFE] \label{cy}\

\ni
Let $\g_0:\R\to \R^3$ be a curve parametrized by arc-length, $g_0$ a $p$-frame along $\g_0$, and  $r_1, r_2$ the corresponding principal curvatures. Let $\phi_0\in SU(2)$ such that $g_0(0)=\Ad(\phi_0)[\d]$.  Suppose $q:\R^2\to \C$ is a solution of the NLS with $q(\cdot,0)=\frac{1}{2}(r_1+ ir_2)$. Let $E$ be the frame of $q$ with $E(0,0,\l)= \phi_0$, and $\a=\frac{\p E}{\p\l}E^{-1}\big|_{\l=0}$. Then $\g(x,t)= \a(x,t)+\g_0(0)$ is the solution of 
\beq\label{ea}
\bca \g_t= \frac{1}{2} \g_x\times \g_{xx}, \\ \g(x, 0)= \g_0(x).\eca
\eeq
\ethm

If the principal curvatures of a curve $\g:\R\to \R^3$ along an orthonormal normal frame $(e_1, e_2)$ are rapidly decaying functions, then so are the principal curvatures along any orthonormal normal frame.  So the property of rapidly decaying principal curvatures is independent of the choice of orthonormal normal frame.  

The Cauchy problem for the NLS with rapidly decaying smooth initial data on $\R$ was solved by the inverse scattering method:

\bthm\label{ci} (\cite{ZS72}, \cite{BeaCoi84}, \cite{BeaCoi85}) Given $q_0\in \calS(\R, \C)$, then  the Cauchy problem of the NLS with initial data $q_0$ has a unique global smooth solution $q:\R^2\to \C$. Moreover,  $q(x, t)$ is rapidly decaying in $x$. 
\ethm

\bcor\label{cz} Let $\g_0:\R\to \R^3$ be a smooth curve with rapidly decaying principal curvatures $k_1, k_2$.  Then \eqref{ea} has a unique smooth solution $\g:\R^2\to \R^3$.
\ecor

\bs
\bs
\section{The VFE on closed curves}\label{du}

In this section, we consider solutions of the VFE, $\g_t=\frac{1}{2}\g_x\times \g_{xx}$, from $S^1\times \R$ to $\R^3$.  We explain the Hasimoto transform for the VFE on closed curves and use the solution of Cauchy problem for the NLS with periodic initial data to solve the periodic Cauchy problem for the VFE. 

The following is a known classical result and we include a proof here. 

\bprop\label{ee} Let $\g:S^1\to \R^3$ be a closed curve parametrized by the arc-length, $\tau$ the torsion, and $R(c)\in SO(2)$ the normal holonomy.   Then  $c=-\oint \tau\rd x$.  
\eprop

\begin{proof}
Let $(e_0, \rn, \rb)$ be the Frenet frame, and $(e_0, e_1, e_2)$ a p-frame along $\g$, and $e_1= \cos \o \rn+ \sin \o \rb$ and $e_2=-\sin \o \rn + \cos \o \rb$. Then $\li (e_1)_x,e_2\ri= \o_x+ \tau$.  
\end{proof}

The following Proposition is also known and the proof is elementary (included).

\bprop\label{cl}
If $\g:S^1\times\R\to \R^3$ is a solution of the VFE, $\g_t= \frac{1}{2} \g_x\times \g_{xx}$, then the normal holonomy of $\g(\cdot, t)$ is independent of $t$.
\eprop

\begin{proof}
Let $f(\cdot,t)=(e_0, \rn, \rb)(\cdot, t)$ be the Frenet frame along $\g(\cdot, t)$, and $A:=f^{-1}f_x= k(e_{21}-e_{12}) + \tau(e_{32}- e_{23})$ and $B= f^{-1} f_t$.  A direct computation implies that $(e_0)_t= \frac{k_x}{2} \rb - \frac{k\tau}{2} \rn$. So 
$$B= \bpm 0& k\tau/2 & -k_x/2\\ -k\tau/2&0& -\xi\\ k_x/2 & \xi &0\epm$$
for some $\xi$.  Equate the $31$-th entry of $A_t= B_x+[A,B]$ to see that $\xi= \frac{k_{xx}}{2k} -\frac{\tau^2}{2}$.  Compare the $32$-th entry to see that $\tau_t= \xi_x+\frac{kk_x}{2}$.  Hence $(\oint \tau(x,t)\rd x )_t=0$. By Proposition \ref{ee}, the normal holonomy for $\g(\cdot, t)$ is independent of $t$. 
\end{proof}

Suppose $\g:S^1\times \R\to \R^3$ is a solution of $\g_t= \frac{1}{2} \g_x\times \g_{xx}$ parametrized by arc-length, and the normal holonomy of $\g(\cdot, 0)$ is trivial. Then the $p$-frame $g(x,t)$ and solution $q(x,t)$ of the NLS constructed in Theorem \ref{bb} is periodic in $x$.  But if the normal holonomy of $\g(\cdot, t)$ is not trivial, then the $g$ and $q$ obtained in Theorem \ref{bb} are not periodic in $x$. Conversely, if $q(x,t)$ is a solution of the NLS that is periodic in $x$, the solution of the VFE constructed in Theorem \ref{bda} may not be closed because there is a period problem when we solve frames for the Lax pair.  To proceed further, we need the following one-parameter Pohlmeyer-Sym construction.

\bthm \label{bc} (\cite{CaIv07}) Suppose $E(x,t,\l)\in SL(2,\C)$ is a frame for the solution $q$ of the NLS.  Given $\l_0\in \R$, let 
\begin{gather}
\a= \frac{\p E}{\p\l} E^{-1}\big|_{\l=\l_0}, \label{aa}\\
\g(x,t)= \a(x-2\l_0 t, t). \label{ab}
\end{gather}
Then $\g$ is a solution of \eqref{cr} on $\R^3$ parametrized by arc-length. 
 \ethm

\begin{proof}
Let $\phi(x,t)= E(x,t,\l_0)$. Then 
$$\bca \phi^{-1}\phi_x= a\l_0+ u,\\ \phi^{-1}\phi_t= a\l_0^2+ u\l_0+Q_{-1}(u),\eca$$
where $a=\diag(i,-i)$, $u=qe_{12}-\bar q e_{21}$, and $Q_{-1}(u)= \frac{i}{2}\bpm -|q|^2 & q_x\\ \bar q_x & |q|^2\epm$.  Use $E^{-1}\rd E= \o_2$ to see that 
$$\bca \a_x= \phi a\phi^{-1}, \\ \a_t= \phi (2a\l_0+ u)\phi^{-1}.\eca$$
So we have $e_0(x,t):= \g_x(x,t)= \phi a\phi^{-1}(x-2\l_0 t, t)$.  Set 
$$e_1(x,t)= -(\phi c\phi^{-1})(x-2\l_0 t, t), \quad e_2(x,t)=(\phi b\phi^{-1})(x-2\l_0 t, t).$$
Let $g=(e_0, e_1, e_2)$. Use a proof similar to Theorem \ref{bda} to get
\beq\label{eta}
 g^{-1}g_x= \bpm 0& -2q_1 &-2q_2\\ 2q_1 &0& -2\l_0\\ 2q_2 & 2\l_0 & 0\epm.
 \eeq
So the mean curvature is $H(\g)= 2q_1 e_1+ 2 q_2 e_2$.  Use the formula for $\a_x$, $\a_t$ and the definition of $\g$ to get 
$$\g_t= -2\l_0 \a_x + \phi(2\l_0+ u)\phi^{-1}= \phi u \phi^{-1}= q_1\phi b\phi^{-1}+ q_2\phi c \phi^{-1}= q_1 e_3- q_2e_1,$$
which is equal to $\frac{1}{2}\g_x\times \g_{xx}$.

Next we compute $g^{-1}g_t$.  Note that $(e_0)_t= (\g_x)_t= (\g_t)_x=( q_1 e_2-q_2 e_1)_x$. Use \eqref{eta} to get
$$(e_0)_t= -((q_2)_x + 2\l_0 q_1)e_1+ ((q_1)_x- 2\l_0 q_2)e_2.$$
Since $(e_2)_t= \phi[-2\l_0 \phi^{-1}\phi_x + \phi^{-1}\phi_t, b]\phi^{-1}$. A direct computation implies that 
$$\li (e_2)_t, e_1\ri= q|^2 + 2\l_0^2.$$ This implies that 
 \beq\label{etb} g^{-1} g_t= \bpm 0 & (q_2)_x+ 2\l_0 q_1 & -(q_1)_x+ 2\l_0 q_2\\ -(q_2)_x-2\l_0 q_1 & 0& |q|^2 + 2\l_0^2\\ (q_1)_x-2\l_0 q_2 & -(|q|^2 + 2\l_0^2) &0\epm.
\eeq
\end{proof}

Given a solution $\g:S^1\times \R\to \R^3$ of the VFE, we have seen that the normal holonomy of $\g(\cdot, t)$ is $R(-2\pi c_0)$ independent of $t$.  We use Example \ref{es}  to construct a $k:S^1\times \R\to SO(3)$ such that $k(\cdot, t)$ is a periodic h-frame along $\g(\cdot, t)$ for each $t$, i.e., $k$ satisfies \eqref{eta} with $\l_0= \frac{1}{2} c_0$.  The proof of the above Theorem tell us how to rotate the h-frame $k(\cdot, t)$ to make the new frame periodic and satisfy \eqref{etb}. This gives a solution of the NLS periodic in $x$.  In fact, we obtain the following.

\bthm\label{fu} Let $\g(x,t)$ be a solution of $\g_t=\frac{1}{2} \g_x\times \g_{xx}$ parametrized by arc-length periodic in $x$ with period $2\pi$, and $R(-2\pi c_0)\in SO(2)$ the normal holonomy for $\g(\cdot, t)$ (by Proposition \ref{cl}, $c_0$ is a constant).  
Then there exists $h:S^1\times \R\to SO(3)$ such that $h(\cdot, t)$ is a periodic $h$-frame for $\g(\cdot, t)$ and $h$ satisfies \eqref{eta} and \eqref{etb} with $\l_0= \frac{c_0}{2}$, i.e., 
\beq\label{ct}
\bca h^{-1}h_x= \bpm 0& -k_1& -k_2\\ k_1& 0& -c_0\\ k_2 & c_0 &0\epm,\\
h^{-1}h_t= \bpm 0 & \frac{(k_2)_x}{2} +\l_0 k_1 & -\frac{(k_1)_x}{2}+\l_0 k_2\\ -\frac{(k_2)_x}{2}- \l_0k_1 &0 & \frac{k_1^2+ k_2^2}{4} +\frac{c_0^2}{2}\\ \frac{(k_1)_x}{2}-\l_0 k_2 &  - \frac{k_1^2+ k_2^2}{4} -\frac{c_0^2}{2}&0\epm.
\eca
\eeq
Moreover,  $q(x,t)=\frac{1}{2}(k_1+ ik_2)(x+c_0t, t)$ is a solution of the NLS periodic in $x$. 
\ethm

\begin{proof}
Since the normal holonomy of $\g(\cdot, t)$ is $R(-2\pi c_0)$, by Example \ref{es}, there exists a periodic h-frame $k$ satisfying 
\beq\label{cta}
A:=k^{-1}k_x= \bpm 0& -r_1 &-r_2\\ r_1 & 0 &-c_0\\ r_2 & c_0 &0\epm.
\eeq
 Write $k=(e_0, e_1, e_2)$. Use the flow equation for $\g$ and \eqref{cta} to see that 
$$(e_0)_t = (\g_x)_t= (\g_t)_x =  \frac{1}{2} (((r_1)_x-c_0r_2)e_2-((r_2)_x+ c_0r_1)e_1).$$ 
So we have 
$$B:= k^{-1}k_t= \bpm 0& \frac{1}{2} (r_2)_x+c_0r_1 & -\frac{1}{2} (r_1)_x+c_0r_2\\ -\frac{1}{2} (r_2)_x-c_0r_1 & 0 &-\xi\\ \frac{1}{2}{(r_1)_x}-c_0r_2& \xi &0\epm$$
for some $\xi$.  Use $A_t= B_x+ [A,B]$ and a direct computation to get $\xi_x=-\frac{1}{4} (r_1^2+r_2^2)_x$. So $\xi= -\frac{1}{4}(r_1^2+ r_2^2) +b(t)$ for some function $b:S^1\to \R$.  Let $\o(t)= -\int_0^t b(s) \rd s$, and $h(x,t)= k(x,t)\diag(1, R(\o(t)-\frac{c_0^2 t}{2}))$. Then $h(\cdot, t)$ is still a periodic h-frame along $\g(\cdot, t)$ and $h$ satisfies \eqref{ct} with 
$$(k_1+ ik_2)(x,t)= e^{i(\frac{c_0^2t}{2}- \o(t))} (r_1+ i r_2)(x,t).$$ A direct computation shows that $q(x,t)= \frac{1}{2}(k_1+ik_2)(x+ c_0t, t)$ is a solution of the NLS periodic in $x$. 
\end{proof}

\bcor  Let $\g:S^1\times \R\to \R^3$ be a solution of $\g_t= \frac{1}{2}\g_x\times \g_{xx}$ parametrized by arc-length with normal holonomy $R(-2\pi c_0)$. Suppose $h:S^1\times \R\to SO(3)$ satisfies the condition that $h(\cdot, t)$ is a periodic h-frame along $\g(\cdot, t)$ for each $t$, i.e., $h^{-1}h_x= \bpm 0& -r_1& -r_2\\ r_1&0&-c_0\\r_2& c_0&0\epm$.  Then there exists $\rho:\R\to \R$ such that $q(x,t)= \frac{1}{2}e^{i\rho(t)} (r_1+ir_2)(x+c_0t, t)$ is a solution of the NLS and $q$ is periodic in $x$. 
\ecor

Let $q_0:S^1\to \C$ be a smooth function.  It follows from results in  \cite{Its76} (also in \cite{Bour93}) that the following periodic Cauchy problem for NLS has a global smooth solution $q:S^1\times \R\to \C$,
$$\bca q_t= \frac{i}{2}(q_{xx}+ 2|q|^2 q),\\ q(x,0)=q_0(x).\eca$$
As a consequence of Theorems \ref{bc} and \ref{fu}, we can use the solution of the periodic Cauchy problem for the NLS to solve the periodic Cauchy problem for the VFE:

\bthm[Periodic Cauchy problem for the VFE]\

\ni Suppose $\g_0: S^1\to su(2)\simeq \R^3$ is a smooth closed curve with normal holonomy $R(-2\pi c_0)$,  $h_0$ is a  periodic $h$-frame along $\g_0$ with$$h_0^{-1}(h_0)_x= \bpm  0 & -r^0_1& -r^0_2\\ r^0_1 &0& -c_0\\ r_2 & c_0 &0\epm.$$ 
Let $\phi_0\in SU(2)$ such that $\Ad(\phi_0)[\d]= h_0(0)$, where $\d=(a,-c,b)$. 
Assume that $q(x,t)$ is the solution of the periodic Cauchy problem for the NLS with $q(\cdot,0)=\frac{1}{2}( r_1^0+ ir_2^0)$.  Let $E(x,t,\l)$ be a frame of the solution $q$ of the NLS such that $E(0,0,\l)= \phi_0$, and $\eta=\frac{\p E}{\p\l} E^{-1}|_{\l= c_0/2}$. Then $\g(x,t)= \eta(x-c_0t, t)+\g_0(0)$ is a solution of $\g_t= \frac{1}{2}\g_x\times \g_{xx}$ with initial data $\g(\cdot,0)= \g_0$ and $\g(x,t)$ is periodic in $x$.
\ethm

\bs
\section{B\"acklund transformations for the VFE\/}\label{dv}

B\"acklund transformations (BT) for the NLS are well-known in the literature, they are either written in terms of a linear system or a system of Ricatti equations. We have seen in section \ref{ca} that the Hasimoto transform is a bijection. So it should be easy to translate the BT theory of the NLS to the curve flow VFE. But the inverse of the Hasimoto transform is given by the Pohlmeyer-Sym formula.  This means that we need to know how the frames are changing when we apply BT. But this was given in \cite{TU00}.  So we can use Theorem \ref{bda} to construct BT for the VFE. We also give an algorithm to construct infinitely many families of explicit solutions of the VFE corresponding to pure soliton solutions of the NLS.  

To use the Hasimoto transform and the BT for the NLS to construct BT for the VFE, we need to have a formula for the frames of the new solution $\ti q$ obtained by applying BT to $q$.  The BT constructed in \cite{TU00} does have this formula. We give a brief review next. 

Given $\a\in \C\backslash \R$,  a Hermitian projection $\pi$ of $\C^2$ (i.e., $\pi^*=\pi$ and $\pi^2=\pi$), let
\beq\label{eg}
f_{\a, \pi}(\l)= \I + \frac{\a-\bar\a}{\l- \a} \pi^\perp,
\eeq
where $\pi^\perp= \I-\pi$.  Then $f_{\a,\pi}$ satisfies the $SU(2)$-reality condition, i.e.,  
$$f_{\a, \pi}^{-1}(\l)= {f_{\a,\pi}(\bar\l)}^\ast.$$

\bthm\label{bg} [Algebraic BT for NLS] (\cite{TU00})\

Let $E(x,t,\l)$ be a frame of a solution $u=\bpm 0& q\\ -\bar q&0\epm$ of the NLS, $\pi$ the Hermitian projection of $\C^2$ onto $\C v$, and $\a\in \C\backslash \R$.  Let $\ti v(x,t)= E(x,t,\a)^{-1}(v)$, and $\ti \pi(x,t)$ the Hermitian projection of $\C^2$ onto $\C \ti v$. Then $\ti u= u+(\bar\a-\a)[\ti \pi, a]$ is a solution  of the NLS. Moreover, $\ti E(x,t,\l)=f_{\a, \pi}(\l) E(x,t,\l) f_{\a, \ti\pi(x,t)}^{-1}(\l)$ is a frame for $\ti u$.
\ethm

Note that the Hermitian projection $\pi$ of $\C^2$ onto $\C \bpm y_1 \\ y_2\epm$ is
$$\pi= \frac{1}{|y_1|^2+ |y_2|^2}\bpm |y_1|^2& y_1\bar y_2\\ \bar y_1 y_2 & |y_2|^2\epm.$$ 
Since $E^{-1}\rd E= \o_2$ and $\ti v(x,t)= E(x,t,\a)^{-1}(v_0)$ in Theorem \ref{bg}, $\ti v$ satisfies the linear system $\rd \ti v= -\o_2(\cdot,\cdot, \a)\ti v$. So  Theorem \ref{bg} gives the well-known BT for the NLS in terms of the linear system for $\ti v$: 

\bthm\label{bga} Let $q$ be a solution of the NLS, $\a\in \C\backslash \R$, and $v_0\in\C^2 \bh 0$. Then the following linear system for $y:\R^2\to \C^2$ has a unique solution with $y(0,0)= v_0$,
\beq\label{da}
 \bca y_x= -(\a a+ u) y,\\ y_t= -(\a^2 a + \a u + Q_{-1}(u))y,\eca
\eeq
where $a=\diag(i, -i)$, $u=\bpm 0& q\\ -\bar q &0\epm$, and $ Q_{-1}(u)=\frac{i}{2} \bpm -|q|^2& q_x\\ \bar q_x & |q|^2\epm$.  
Moreover, if $y=(y_1, y_2)^t$ is a solution of \eqref{da}, then
\beq\label{bgb}
\ti q= q + 2i(\a-\bar\a) \frac{y_1 \bar y_2}{||y||^2}
\eeq
is again a solution of the NLS. 
\ethm 

System \eqref{da} explicitly written down is
\beq\label{daa}
\bca
y_x =\bpm -i\a &- q\\ \bar q & i\a\epm y,\\
y_t= \bpm -i\a^2+\frac{i}{2} |q|^2& -\a q -\frac{i}{2} q_x\\ \a \bar q -\frac{i}{2} \bar q_x& i\a^2 -\frac{i}{2} |q|^2 \epm y
\eca
\eeq

Next we review the classical relation between Ricatti equations and systems of linear first order equations. A direct computation shows that if 
\beq\label{cka}
\bpm y_1\\ y_2\epm_x= \bpm A& B\\ C& -A\epm \bpm y_1\\ y_2\epm,
\eeq
then $p=\frac{y_1}{y_2}$ is a solution of 
\beq\label{ckb}
p_x= -Cp^2+ 2Ap + B.
\eeq
Conversely, if $p$ is a solution of \eqref{ckb} and $z$ is a solution of $z_x= Cpz- Az$ then $y=(pz, z)^t$ is a solution of \eqref{cka}.  Combine these with Theorem \ref{bga} we obtain the following well-known BT as a system of compatible Ricatti equations for $p$.

\bthm [Analytic BT for NLS] \label{bh}\

\ni Let $q$ be a solution of the NLS, and $u=\bpm 0& q\\ -\bar q&0\epm$, and $\a\in \C\bh \R$ a constant. Then the following system of differential equations is solvable for $p(x,t)$:
\beq\label{bt}
{\rm(BT)\/}_{q,\a} \,\,\bca p_x= -\bar q p^2 - 2i\a p - q,\\ 
p_t= (\frac{i}{2} \bar q_x -\a \bar q) p^2 + i(|q|^2 -2\a^2) p -(\a q +\frac{i}{2} q_x).\eca 
\eeq
Moreover, if $p$ is a solution of \eqref{bt}, then $\ti q= q+2i(\a-\bar \a) \frac{p}{1+|p|^2}$ is again a solution of the NLS.  
\ethm

In fact, BT comes from the action of the group of rational loops.  First we recall the generator theorem of Uhlenbeck.

\bthm (\cite{Uhl89})
The group $L^r(GL(2))$ of rational maps from $\C\cup \{\infty\}$ to $GL(2,\C)$ satisfying the $SU(2)$-reality condition $f(\bar\l)^*f(\l)=\I_2$ and $f(\infty)=\I_2$ is generated by 
$$\{f_{\a, \pi}\n \a\in \C\bh \R, \pi \,\, {\rm a\, Hermitian\, projection\, of \,\,} \C^2\}.$$
\ethm

Given a solution $q$ of the NLS, we call the solution of 
$$E^{-1}E_x= a\l +u, \quad E^{-1}E_t= a\l^2+ u\l + Q_{-1}(u), \quad E(0,0,\l)=\I,$$
the {\it normalized frame of $q$\/}.

\bthm (\cite{TU00}) Let $E(x,t,\l)$ be the normalized frame of a solution $u=\bpm 0& q\\ -\bar q&0\epm$ of the NLS, and $f\in L^r(GL(2))$. Then there exist unique $\ti E(x,t, \l)$ and $\ti f(x,t,\l)$ satisfying 
$$f(\l) E(x,t,\l)= \ti E(x,t,\l) \ti f(x,t,\l),$$
such that $\ti E(x,t,\l)$ is holomorphic for $\l\in \C$ and $\ti f(x,t, \cdot)\in L^r(GL(2))$. Moreover, 
\ben
\item $\ti E^{-1}\ti E_x= a\l + \ti u$, where  $\ti u=\ti q e_{12}- \overline{\ti q} e_{21}$ and $\ti q$ is a solution of the NLS,
\item $f\ast q:= \ti q$ defines an action of  $L^r(GL(2))$ on the space of solutions of the NLS,
\item the solution $\ti q$ constructed in Theorem \ref{bg} from the normalized frame $E$ of $q$ is equal to $f_{\a,\pi}\ast q$,
\item if $q:\R^2\to \C$ is a smooth solution of the NLS such that $q(x,t)$ is rapidly decaying in $x$, then $f_{\a,\pi}\ast q$ is a smooth solution of the NLS on $\R^2$ that is rapidly decaying in $x$.
\een
\ethm

The Permutability Theorem for BT is a consequence of the following relation of generators of $L^r(GL(2))$.

\bprop (\cite{TU00}) Given $\a_1,\a_2\in \C$ such that $\a_1\not=\a_2, \bar\a_2$ and $\pi_i$  Hermitian projections of $\C^2$ onto $V_i$ for $i=1,2$.  Let
$$W_1= f_{\a_2,\pi_2}(\a_1)(V_1), \quad W_2=f_{\a_1, \pi_1}(\a_2)(V_2),$$
and $\tau_i$ the Hermitian projection of $\C^2$ onto $W_i$ for $i=1,2$. 
Then 
$$f_{\a_1,\tau_1} f_{\a_2,\pi_2}= f_{\a_2, \tau_2}f_{\a_1, \pi_1}.$$
\eprop

\bthm[Permutability for BT]\label{fv} (\cite{TU00})\

\ni Let $\a_1, \a_2\in \C\backslash \R$ such that $\a_1\not= \a_2, \bar \a_2$, $q$ a solution of the NLS,  $p_1$ and $p_2$ solutions of \eqref{bt} with parameter $\a_1, \a_2$ respectively. Let $q_i$ be the solution of the NLS constructed from $p_i$, and $\ti \pi_i(x,t)$ the Hermitian projection onto $\C (p_i(x,t),1)^t$ for $i=1,2$. Let $\ti\tau_1(x,t)$ and $\ti \tau_2(x,t)$ be the Hermitian projections of $\C^2$ onto $\C\xi(x,t)$ and $\C \eta(x,t)$ respectively, where
\begin{align*}
&\xi(x,t)=\bpm \xi_1(x,t)\\ \xi_2(x,t)\epm=f_{\a_2,\ti\pi_2(x,t)}(\a_1)\bpm p_1(x,t)\\ 1\epm, \\
&\eta(x,t)= \bpm \eta_1(x,t)\\ \eta_2(x,t)\epm= f_{\a_1, \ti\pi_1(x,t)}(\a_2)\bpm p_2(x,t)\\ 1\epm.
\end{align*}
Then $\ti p_1= \xi_1/\xi_2$ is a solution of (BT)$_{q_2, \a_1}$, $\ti p_2=\eta_1/\eta_2$ is a solution of (BT)$_{q_1, \a_2}$, and 
$$q_3= q_2 +2i(\a_1-\bar\a_1)\frac{\ti p_1}{1+|\ti p_1|^2}= q_1 + 2i(\a_2-\bar\a_2)\frac{\ti p_2}{1+|\ti p_2|^2}$$
is a solution of the NLS.
\ethm

A direct computation implies that 
$$\ti p_1= \frac{(\a_1-\bar \a_2)p_1+ (\a_1-\a_2) p_1 |p_2|^2 -(\a_2-\bar \a_2) p_2}{(\a_1-\a_2) +(\a_1-\bar\a_2)|p_2|^2-(\a_2-\bar \a_2) p_1\bar p_2},$$
and $\ti p_2$ is given by the same formula but interchange the indices $1$ and $2$.

  As a consequence of Theorems \ref{bb},  \ref{bda}, \ref{bd} and \ref{bg}, we obtain the following.

\bthm\label{dm}
Let $\g(x,t)$ be a solution of $\g_t=\frac{1}{2} \g_x\times \g_{xx}$ parametrized by arc-length, and $g, q, \phi, E$ as in Theorem \ref{bd}.  Let $\a\in \C\bh \R$,  $v_0\in\C^2\bh 0$, and $\ti \pi(x,t)$ the projection of $\C^2$ onto $\C y(x,t)$, where $y=\bpm y_1\\ y_2\epm= E(\cdot, \cdot,\a)^{-1}(v_0)$. Let $\phi(x,t)= E(x,t,0)$.  Then 
$$\ti \g= \g-\frac{(\a-\bar\a)}{|\a|^2}\, \phi\ti \pi_\ast \phi^{-1}$$
is a new solution of the $\ast$-MCF on $\R^3$, where $\ti\pi_\ast= \ti\pi-\frac{1}{2}\I_2$ is the traceless part of $\ti \pi$. 
Moreover, let $\ti\phi= \phi (\ti \pi +\frac{\a}{\bar\a} \ti\pi^\perp)$, and $\ti g=\Ad(\ti\phi)[\d]$, then we have
\ben
\item[(i)] $\ti g(\cdot, t)$ is a p-frame along $\ti \g(\cdot, t)$ with principal curvatures $\ti k_1, \ti k_2$,
\item[(ii)] $\ti q= \frac{1}{2} (\ti k_1+ i\ti k_2)= q+ 2i(\a-\bar \a) \frac{y_1\bar y_2}{|y_1|^2+|y_2|^2}$ is a solution of the NLS,
\item[(iii)] $\ti E(x,t,\l) = E(x,t,\l) f^{-1}_{\a, \ti\pi(x,t)}(\l)$ is a frame for $\ti \g$.
\een
\ethm

Next use the p-frame $g=(e_0, e_1, e_2)= \Ad(\phi)[\d]$ and the solution $y$ of the linear system \eqref{da} to write the BT for the VFE.

\bthm\label{cj} Let $\g$ be a solution of $\g_t= \frac{1}{2}\g_x\times \g_{xx}$ parametrized by the arc-length parameter, $g=(e_0, e_1, e_2)$, $k_i(\cdot, t)$ the principal curvature along $e_i(\cdot, t)$ for $i=1,2$,  and $q=\frac{1}{2}(k_1+ik_2)$ the solution of NLS as in Theorem \ref{bb}. Let $\a\in \C\bh \R$, and $y=(y_1, y_2)^t$ a solution of \eqref{da}. Then 
\beq\label{dl}
\ti \g= \g + \frac{\Im (\a)}{||\a y||^2}\left( (|y_2|^2- |y_1|^2) e_0 + (y_1 \bar y_2+ \bar y_1 y_2) e_1 + \frac{1}{i}(y_1 \bar y_2 - \bar y_1 y_2) e_2\right)
\eeq
is a solution of $\g_t= \frac{1}{2}\g_x\times \g_{xx}$.  
\ethm

BT for the VFE can also be written in terms of a system \eqref{bt} for $p$: Let $\g$ be a solution of the VFE and $p$ a solution of \eqref{bt}, and $g=(e_0, e_1, e_2)$, $q$ solution of the NLS as in Theorem \ref{cj}.   If $p$ is a solution of (BT)$_{q,\a}$, then  
\beq\label{dla}
\ti\g= \g+\frac{\Im\a}{|\a|^2(1+|p|^2)} \left((1-|p|^2) e_0+ (p+\bar p) e_1 + \frac{1}{i} (p-\bar p) e_2\right)
\eeq
is a solution of the VFE. 

Similarly, a permutability formula for BT of the $\ast$-MCF on $\R^3$ can be written down as follows:

\bthm[Permutability for the VFE]\label{fx}\

\ni Let $\g, g=(e_0, e_1, e_2)$ and $q$ be as in Theorem \ref{cj}. Let $\a_1, \a_2\in \C\bh 0$ such that $\a_1\not= \a_2, \bar \a_2$, $p_i$ a solution of (BT)$_{q, \a_i}$, and  $\ti \pi_i$ the Hermitian projection onto $\C (p_i, 1)^t$ for $i=1,2$.  Let $\phi:\R^2\to SU(2)$ such that $g= \Ad(\phi)[\d]$, $\phi_i=\phi (\ti\pi_i + \frac{\a_i}{\bar\a_i} \ti \pi_i^\perp)$, and 
$$g_i=(e_0^i, e_1^i, e_2^i)= \Ad(\phi_i)[\d]$$
for $i=1,2$. Let   $\g_i$ be the solution given by $p_i$ as in \eqref{dla}, i.e.,
$$\g_i = \g+\frac{\Im(\a_i)}{|\a_i|^2 (1+|p_i|^2)}((1-|p_i|^2) e_0+ 2\Re(p_i) e_1 + 2\Im(p_i) e_2),$$
for $i=1,2$.  
Let $\ti p_1, \ti p_2$ be as in Theorem \ref{fv}. Then 
\begin{align*}
\g_{12}&= \g_1 +\frac{\Im(\a_2)}{|\a_2|^2(1+|\ti p_2|^2)} ((1-|\ti p_2|^2)e^1_0 + 2\Re(\ti p_2) e^1_1 + 2\Im(\ti p_2) e^1_2)\\
&= \g_2 + +\frac{\Im(\a_1)}{|\a_1|^2(1+|\ti p_1|^2} ((1-|\ti p_1|^2)e^2_0 + 2\Re(\ti p_1) e^2_1 + 2\Im(\ti p_1) e^2_2)
\end{align*}
is a solution of the VFE. 
\ethm 

\brem[Explicit solutions of the VFE]\ 

\ni
Note that the straight line $\g_0(x,t)= ax$ is a stationary solution of the VFE, $\g_t= \frac{1}{2}\g_x\times \g_{xx}$. The corresponding solution of the NLS is  the trivial solution $q=0$. Since $E(x,t,\l)= \exp(a(\l x+ \l^2 t))$ is the frame for the trivial solution $q=0$ of the NLS,  we can apply Theorem \ref{dm} to the stationary solution $\g_0$ to obtain explicit solutions of the VFE corresponding to 1-soliton solutions of the NLS. We can then either apply BT Theorem \ref{dm}  repeatedly $k$ times or use the Permutability Theorem \ref{fx} for the VFE to obtain explicit solutions of the VFE  that correspond to the explicit pure $k$-soliton solutions of the  NLS. 
\erem

\bs
\section{The Hamiltonian aspects of the VFE}\label{ef}

Let $\calC=\{\g\in C^\infty(\R, \R^3)\n ||\g_x||=1\}$, $\fR(3)$ the group of rigid motions of $\R^3$, and $V_1= \{y e_{12}-\bar y e_{21}\n y\in \C\}$.   Let 
$$\fy:\calC/\fR(3)\to C^\infty(\R, V_1)/S^1$$ be the map defined by $\fy([\g])= [u]$, where $u= q e_{12}- \bar q e_{21}$ is constructed as follows: Choose a p-frame $g$ along $\g=(e_0, e_1, e_2)$, let $k_i$ be principal curvature along $e_i$ for $i=1,2$. Then $q=\frac{1}{2}(k_1 + ik_2)$. Here $[\g]$ and $[u]$ are the $\fR(3)$-orbit of $\g$ and the $S^1$-orbit of $[u]$ respectively.   It follows from  the existence and uniqueness of ordinary differential equations that $\fy$ is a bijection.  Let
$\calC_s$ denote the subset of all $\g\in \calC$ that have rapidly decaying principal curvatures,   and $\calS(\R, V_1)$ the space of rapidly decaying maps from $\R$ to $V_1$. Then $\fy$ maps $\calC_s/ \fR(3)$ bijectively onto $\calS(\R, V_1)/S^1$.  The formulas \eqref{bn1} and \eqref{bn2} define two Poisson structures on $\calS(\R, V_1)$. By Proposition \ref{dd}, $\{F_1, F_2\}_i$ is $S^1$-invariant  if $F_1, F_2$ are $S^1$-invariant functionals on $\calS(\R, V_1)$. So these two Poisson structures are in fact defined on $\calS^\infty(\R, V_1)/S^1$ or on $S^1$-invariant functionals on $\calS(\R, V_1)$.   Since $\fy$ is a bijection, the map $F\mapsto \hat F= F\circ \fy$ gives a bijection between the space of $S^1$-invariant functionals on $\calS(\R, \C)$ and the $\fR(3)$-invariant functionals on $\calC_s$.   So the pull back  of the Poisson structures $\{\, ,\}_i$ on $\calS(\R, \C)$ by $\fy$  gives a $\fR(3)$-invariant Poisson structure on $\calC_s$:
$$\{F\circ\fy, H\circ \fy\}_i^\wedge(\g)= \{F, H\}(\fy(\g)), \quad i=1,2.$$
  We have seen in section \ref{ca} that the VFE is invariant under the rigid motion group $\fR(3)$, the NLS is invariant under the action of $S^1$, and the VFE flow corresponds to the NLS flow under $\fy$. So we can use $\fy$ to translate the Hamiltonian aspects of the NLS to get the Hamiltonian aspects of the VFE.  
In particular, we write down the commuting conservation laws and the higher commuting curve flows for the VFE.   

Let $H_j(u)=-\frac{1}{j}\oint \li Q_{-j}(u),a\ri \rd x$ be as defined by \eqref{bq}, and $\hat H_j:\calC\to \R$ the functional on $\calC_s$ defined by $\hat H_j(\g)= H_j(u)$, where $\fy([\g])= [u]$, $u=qe_{12}-\bar q e_{21}$, $q=\frac{1}{2}(k_1+ ik_1)$, and $k_1, k_2$ are principal curvatures defined by a p-frame along $\g$. For example,
\begin{align*}
\hat H_1(\g)&= \frac{1}{2} \oint |q|^2\rd x= \frac{1}{8} \oint k^2\rd x,\\
\hat H_2(\g)&= -\frac{1}{4} \oint \Im(\bar q q_x)\rd x = \frac{1}{16} \oint k^2\tau \rd x,
\end{align*}
where $q=\frac{1}{2}(k_1+ik_2)$, $k_1, k_2$ the principal curvatures with respect to a p-frame, and $k, \tau$ the curvature and torsion of $\g$.

Next we use the map $\fy$ to write down the curve flow that correspond to the $j$-th flow in the $SU(2)$-hierarchy. Write $q=\frac{1}{2}(q_1+iq_2)$ and
$$Q_{-j}(u)=A_{-j} (q) a+ B_{-j}(q)b + C_{-j}(q) c,$$ where  $u=q e_{12}- \bar q e_{21}$.  For example, write $q=q_1+iq_2$, then 
\begin{align*}
&Q_0(u)=q_1 b+ q_2 c,\\
&Q_{-1}(u)= -\frac{1}{2} |q|^2 a-\frac{1}{2} (q_2)_x b + \frac{1}{2} (q_1)_x c,\\
&Q_{-2}(u)= \frac{1}{2} \Im(\bar q q_x) a -\frac{1}{4} ((q_1)_{xx} + 2|q|^2 q_1) b -\frac{1}{4} ((q_2)_{xx} + 2|q|^2 q_2) c.
\end{align*}  
For $j\geq 1$, we call the following curve flow on $\R^3$ {\it the $j$-th curves flow\/},
\beq\label{de}
\g_{t_j}= A_{-(j-2)}(q) e_0 -C_{-(j-2)}(q) e_1 + B_{-(j-2)}(q) e_2,
\eeq
where $q=\frac{1}{2}(k_1+ ik_2)$ and  $k_1, k_2$ are the principal curvatures with respect to a p-frame $(e_0, e_1, e_2)$ along $\g$.  Next we write the first four flows for $1\leq j\leq 4$ invariantly: 
\begin{align}
&\g_t= \g_x, \\
& \g_t=\frac{1}{2}(e_0\times H(\g)), \\
&\g_t= -\frac{1}{8} ||H(\g)||^2e_0 - \frac{1}{4} \K_{e_0}^\perp H(\g),  \label{eb1}\\
&\g_t=\frac{1}{8} k^2\tau \, e_0- \frac{1}{8} e_0\times (\K_{e_0}^\perp)^2 H(\g) - \frac{||H(\g)||^2}{16} e_0\times  H(\g),\label{eb2}
\end{align}
where $e_0$ is the unit tangent to the curve and $k, \tau$ are the curvature and torsion.  Here we use the identity $\Im(\bar q q_x)= -\frac{1}{4} k^2\tau$. Also note that $e_0\times \xi= \ast_\g(\xi)$ for any normal field $\xi$ along $\g$.  
The following proposition states that the right hand side of \eqref{de} is independent of the choice of p-frames for all $j\geq 1$:

\bprop The flow \eqref{de} is a geometric curve flow, i.e., the right hand side of \eqref{de} is independent of the choice of the p-frames.
\eprop

\begin{proof} if both $(e_0, e_1, e_2)$ and $(e_0, \ti e_1, \ti e_2)$ are p-frames along $\g$, then there is a constant $c\in \R$ such that $(\ti e_1, \ti e_2)= (e_1, e_2)R(c)$.  Let $k_i$ and $\ti k_i$ be the principal curvatures of $\g$ along $e_i$ and $\ti e_i$ respectively. Then $\ti k_1+ i\ti k_2= e^{-ic} (k_1+ ik_2)$.  So $\ti q= e^{-ic} q$.  Let $u=qe_{12} -\bar q e_{21}$, and $\ti u= \ti q q_{12}- \overline{\ti q} e_{21}$. Then $\ti u= e^{-ca/2} u e^{c a/2}$.  By \eqref{dj}, $Q_{-j}(\ti u)= e^{-ca/2} Q_{-j}(u) e^{ca/2}$. This implies that $A_j(\ti q)= A_j(q)$ and 
$$-C_j(\ti q) \ti e_1 + B_j(\ti q) \ti e_2= -C_j(q) e_1 + B_j(q) e_2.$$
So the right hand side of \eqref{de} is independent of the choice of p-frames. 
\end{proof}

\bprop  If a geometric curve flow written in terms of a p-frame $(e_0, e_1, e_2)$ is $\g_t= \sum_{i=0}^2 \xi_i(k_1, k_2) e_i$, then the curve flow preserves the arc-length parameter if and only if  $e_0(\xi_0)= k_1\xi_1 + k_2\xi_2$.
\eprop

\begin{proof}  Let $s$ denote the arc-length parameter, then $\xi_x= ||\g_x|| \xi_s$. 
We compute directly to get
\begin{align*}
&\frac{1}{2}(\li \g_x, \g_x\ri)_t= \li \g_{xt}, \g_x\ri =\li (\g_t)_x, \g_x\ri \\
&\quad =||\g_x||^2 \li ( \xi_0 e_0+ \xi_1 e_1+ \xi_2 e_2)_s, e_0\ri = ||\g_x||^2 ((\xi_0)_s - k_1 \xi_1 -k_2\xi_2).
\end{align*}
This proves the claim.
\end{proof} 

\bcor
The flow \eqref{de} preserves the arc-length parameter.
\ecor

\begin{proof}
Equate the coefficient of $a$ in \eqref{df}, 
$$(Q_{-j}(u))_x + [u, Q_{-j}(u)]=[Q_{-(j+1)}(u), a],$$
 to get $(A_{-j})_s + C_{-j} k_1 -B_{-j} k_2=0$. Since $\xi_0= A_{-j}$, $\xi_1= -C_{-j}$, and $\xi_2= B_{-j}$, we get $(\xi_0)_s= k_1 \xi_1+ k_2 \xi_2$. 
\end{proof}

Use same proofs as in section \ref{ca} to get the following results for the $j$-th curve flow.  

\ben
\item $\fy$ maps bijectively solutions of the $j$-th curve flow \eqref{de} modulo $\fR(3)$ to solutions of the $j$-th flow in the $SU(2)$-hierarchy modulo the $S^1$-action. 
\item The $j$-th curve flow \eqref{de} is the Hamiltonian equation for $\hat H_j$ on $\calC_s$ with respect to $\{\, ,\}^\wedge_2$, and it is also the Hamiltonian equation for $\hat H_{j+1}$ with respect to $\{\,, \}^\wedge_1$.
\item These curve flows commute. In particular, the $j$-th flow commutes with the $\ast$-MCF on $\R^3$. 
\item B\"acklund theory for the $j$-th curve flow can be worked out in a similar manner. For example, Theorem \ref{cj} holds for the $j$-th curve flow if we choose $y$ to be a solution of the linear system $\rd y= -\o_j(x,t,\a) y$, where  
$$\o_j= (a\l+u)\rd x+ (a\l^j+ u\l^{j-1} + Q_{-1}(u)\l^{j-2} + \cdots + Q_{-(j-1)}(u))\rd t$$
 is the Lax pair \eqref{dc} of the $j$-th flow in the $SU(2)$-hierarchy, where $a=\diag(i,-i)$ and $u=qe_{12}-\bar q e_{21}$.
 \een

 Note that changing the tangential component of the curve flow amounts to changing parametrization of the curves.   
Let 
$$X_j(\g)=  A_{-(j-2)}(q) e_0 -C_{-(j-2)}(q) e_1 + B_{-(j-2)}(q) e_2,$$
and  $Y_j(\g)= (X_j(\g))^\perp$, the normal component of $X_j(\g)$, i.e., 
$$Y_j(\g)= -C_{-(j-2)}(q) e_1 + B_{-(j-2)}(q) e_2.$$
 So the curve flow
\beq\label{dg}
\g_t= Y_j(q)= -C_{-(j-2)}(q) e_1 + B_{-(j-2)}(q) e_2
\eeq
geometrically is the same curve flow \eqref{de}.  In other words, the curve at time $t$ of a solution of \eqref{dg} is the same curve at time $t$ of \eqref{de} but with different parametrization. 
For example, the curve flow \eqref{dg} for $1\leq j\leq 4$ are 
\begin{align}
&\g_{t_1}= \g_x, \label{dh1}\\
&\g_{t_2}= \frac{1}{2}\ast_\g H(\g), \label{dh2}\\
&\g_{t_3}= -\frac{1}{4} \K^\perp_{e_0} H(\g), \label{dh3}\\
& \g_{t_4}= -\frac{1}{16} \ast_\g(||H(\g)||^2 H(\g) +2 (\K^\perp_{e_0})^2 H(\g)). \label{dh4}
\end{align}
Note that \eqref{dh2} is the $\ast$-MCF and \eqref{dh3} is the geometric Airy flow.  These flows preserves total arc-length, so we can reparametrize the curves by the arc-length and get the curve flow \eqref{de}. In particular, we get the following.

\bthm The geometric Airy flow \eqref{dh3} preserves total arc-length and the normalized geometric Airy flow \eqref{eb1} preserves the arc-length parameter. Moreover, the normalized Airy flow commutes with the VFE.  
\ethm

\brem 
It is known that the odd flows in the $SU(2)$-hierarchy leaves the condition $\bar q= q$ invariant and the resulting hierarchy is the hierarchy containing the mKdV for $q:\R^2\to \R$, $q_t= -\frac{1}{4}(q_{xxx} + 6q^2 q_x)$.  The condition that $q$ is real means that the torsion $\tau$ of $\g$ is zero if $\fy([\g])=[q]$. So the $(2j+1)$-th curve flow on $\R^3$ leaves plane curves invariant.  For example, the third flow \eqref{de} on $\R^3$ leaves the plan curves invariant and the curve flow on the plane is 
$$\g_t= -\frac{1}{8} k^2 e_0 -\frac{1}{4} e_0(k) e_1,$$ where $k$ is the curvature of the plane curve. This is the normalized curve flow for the geometric Airy flow on $\R^2$, $\g_t= -\frac{1}{4} \K^\perp_{e_0} H(\g(\cdot, t))$. 
\erem

\bs

\section{The time-like $\ast$-MCF on $\R^{2,1}$}\label{dw}

In this section, we explain the relation between the time-like $\ast$-MCF on $\R^{2,1}$ and the defocusing NLS and give B\"acklund transformations for this curve flow.

Recall that  
\begin{align*}
SU(1,1)&=\{g\in SL(2,\C)\n g^*Jg=J\}, \quad {\rm where\,\,} J=\diag(1,-1),\\
su(1,1) &=\left\{\bpm ix & z\\ \bar z & -ix\epm\n x\in \R, z\in \C\right\}.
\end{align*}
For time-like solution of the $\ast$-MCF on $\R^{2,1}$ we identify $\R^{2,1}$ as $su(1,1)$ as follows: Let 
$$\li X, Y\ri = \frac{1}{2} \tr(XY)$$ be the bi-linear on $su(1,1)$. Then 
$$\li X, X\ri = -x^2+|z|^2, \quad {\rm for\,\,} X=\bpm ix & z\\ \bar z& -ix\epm.$$
 Let
\beq\label{cv}
\d_1= (a, b, c), \quad a=\diag(i, -i), \quad b= \bpm 0&1\\ 1&0\epm, \quad c= \bpm 0& i\\-i &0\epm.
\eeq
Then $\d_1$ is an ordered basis of $su(1,1)$ satisfying 
$$-\li a, a\ri =\li b, b\ri =\li c, c\ri =1, \quad \li a, b\ri=\li a, c\ri = \li b, c\ri=0.$$
 We identify $su(1,1)$ as $\R^{1,2}$ by 
$$ax+by+cz\mapsto (x,y,z)^t.$$
Let 
$$D_2=\diag(-1, 1, 1), \quad SO(1,2)=\{A\in GL(3)\n A^t D_2A= D_2\}.$$ Then
$$so(1,2)=\left\{ \bpm 0& y_1 & y_2\\ y_1& 0 & -y_3\\ y_2 & y_3 &0\epm\,\bigg|\, y_1, y_2, y_3\in \R\right\}.$$
Let $\Ad: SU(1,1)\to O(su(1,1))=SO(1,2)$ be the Adjoint representation of $SU(1,1)$ on $su(1,1)$, and 
\beq\label{cva}
\Ad(g)[\d_1]=(gag^{-1}, gbg^{-1}, gcg^{-1})\in SO(1,2).
\eeq 

The following theorem can be proved the same way as Theorem \ref{bb}.

\bthm\label{fm} If $\g:\R^2\to \R^{2,1}$ is a solution of the $\ast$-MCF on $\R^{2,1}$
such that $\li \g_x, \g_x\ri=-1$ (so $\g$ is a time-like solution),
then there exists $g=(e_0, e_1, e_2):\R^2\to SO(1,2)$ such that
$$g^{-1}g_x= \bpm 0& k_1 & k_2\\ k_1&0&0\\ k_2&0&0\epm$$
and $q=\frac{1}{2}(k_1+ ik_2)$ is a solution of the defocusing NLS, $q_t= \frac{i}{2}(q_{xx} - 2|q|^2 q)$.  Moreover, if $\ti g:\R^2\to SO(1,2)$ is another map satisfies these conditions, then there exists a constant $R(c)\in SO(2)$ such that $\ti g= g\diag(1, R(c))$ and the corresponding $\ti k_1, \ti k_2$ satisfies $\ti k_1+\ti k_2= e^{-ic} (k_1+ ik_2)$.
\ethm

\brem Let 
\begin{align*}
&\fR_1(3)= \Iso(\R^{2,1})= SO(1,2)\ltimes \R^{2,1},\\
&\calC_1=\{\g\in C^\infty(\R, \R^{2,1})\n \li \g_x, \g_x\ri =-1\},\\
&V_2=\{ye_{12}+ \bar y e_{21}\n y\in \C\}.
\end{align*} 
Theorem \ref{fm} associate to each solution $\g$ of the time-like $\ast$-MCF an $S^1$-orbit of solutions of the defocusing NLS. Moreover, if $\g_2=\phi\circ g_1$ for some $\phi\in \fR_1(3)$, then the $S^1$-orbits of solutions of the defocusing NLS associated to $\g_1$ and $\g_2$ are the same. Let $\Psi_1$ denote the map from the orbit space of the space of solutions of the time-like $\ast$-MCF on $\R^{2,1}$ under the action of $\fR_1(3)$ and the orbit space of the space of solutions of the defocusing NLS under the action of  $S^1$ defined by $\Psi([\g])=[u]$, where $u=qe_{12}+ \bar q e_{21}$ is a solution of the defocusing NLS associated to $\g$ obtained in Theorem \ref{fm}. 
Then $\Psi_1$ is a well-defined map and is one to one. The same proof of \ref{bda} implies that if $q$ is a solution of the defocusing NLS and $E(x,t,\l)$ the normalized frame for $u=qe_{12}+ \bar q e_{21}$, then $\a=\frac{\p E}{\p\l}E^{-1}|_{\l=0}$ is a solution of the time-like $\ast$-MCF on $su(1,1)\simeq \R^{2,1}$ and  $\Psi_1([\a])= [u]$. This shows that $\Psi_1$ is a bijection and is the analogous Hasimoto transform for the time-like $\ast$-MCF on $\R^{2,1}$.
\erem

To construct BT for the time-like $\ast$-MCF on $\R^{2,1}$, we need to construct 
B\"acklund transformations for the $SU(1,1)$-hierarchy. First we write down the projection of $\C^2$ with respect to the bi-linear form $\li\, ,\ri_1$ on $\C^2$ defined by
$$\li X, Y\ri_1= X^\ast J Y, \quad J=\diag(1,-1).$$  The Adjont of a linear operator $A:\C^2\to \C^2$ is 
$$A^\sharp = JA^* J.$$
(Recall that the Adjoint of $A$ is defined by $\li AX, Y\ri_1= \li X, A^\sharp Y\ri_1$ for all $X, Y\in \C^2$).
Let $\pi\in GL(2,\C)$ satisfying 
\beq\label{eh}
\pi^\sharp=\pi, \quad \pi^2=\pi.
\eeq 
Let $\pi^\perp= \I-\pi$. 

If $\pi$ satisfies \eqref{eh} and $v$ such that $\Im(\pi)=\C v$ and $\li v, v\ri_1\not=0$, then 
\beq\label{fy}
\pi= \frac{1}{\li v, v\ri_1} \, vv^\ast J.
\eeq
We call such $\pi$ the projection onto $\C v$ (note that \eqref{fy} only works when $\li v, v\ri_1\not=0$).  

Let $f_{\a,\pi}(\l)=\I+\frac{\a-\bar\a}{\l-\a} \pi^\perp$ be defined by the same formula \eqref{eg} with $\pi$ satisfying \eqref{eh}.  A direct computation implies that $f_{\a, \pi}$ satisfies the $SU(1,1)$-reality condition, 
$$Jf(\bar\l)^* J= f(\l)^{-1}.$$  BT for the defocusing NLS can be proved in a similar way as the focusing NLS.  We summarize the result in the next theorem. 

\bthm Let $q$ be a solution of the defocusing NLS, $q_t= \frac{i}{2}(q_{xx}- 2|q|^2 q)$, and $E(x,t,\l)$  a frame of $q$. Let $\a\in \C\bh \R$,  $v_0\in \C^2\bh 0$ with $\li v_0,v_0\ri_1\not=0$, and $\pi= \frac{1}{\li v_0, v_0\ri_1} v_0 v_0^* J$. If 
$$y(x,t)=(y_1(x,t), y_2(x,t))^t:= E(x,t,\a)^{-1}(\a)$$ is not a null vector for all $(x,t)$ in an open subset $\calO_0$ of $\R^2$, Then we have the following. 
\ben
\item[(i)] $\ti q= q- 2i(\a-\bar\a) \frac{y_1\bar y_2}{|y_1|^2- |y_2|^2}$ is a solution of the defocusing NLS defined on $\calO_0$.
\item[(ii)] $\ti E(x,t,\l)= f_{\a, \pi}(\l) E(x,t,\l)f_{\a, \ti \pi(x,t)}(\l)^{-1}$ is a frame for $\ti q$, where $\ti \pi=\frac{1}{|y_1|^2-|y_2|^2} y y^*J$. 
\item[(iii)] $y$ satisfies $\rd y= -\o_2(\cdot, \cdot, \a) y$, where 
$$\o_2= (a\l + u) \rd x + (a\l^2+ u\l + Q_{-1}(u))\rd t,$$
$ a= \diag(i,-i)$, $u=\bpm 0& q\\ \bar q&0\epm$, and $Q_{-1}(u)=\frac{i}{2}\bpm |q|^2 & q_x\\ -\bar q_x & -|q|^2\epm$; i.e., 
\beq\label{ev}
\bca y_x= -\bpm i\a& q\\ \bar q& -i\a\epm y,\\
y_t= -\bpm i\a^2+\frac{i}{2} |q|^2 & \a q +\frac{i}{2} q_x\\ \a \bar q -\frac{i}{2} \bar q_x & -i\a^2 -\frac{i}{2} |q|^2 \epm y.
\eca
\eeq
Moreover, if $y$ is a solution of \eqref{ev}, then the formula of $\ti q$ given in (i) is a solution of the defocusing NLS. 
\item[(iv)] The following system for $p$ is solvable,
\beq\label{ej}\bca 
p_x= \bar q p^2- 2i\a p-q,\\ p_t= (\a \bar q- \frac{i}{2} \bar q_x)p^2- i(2\a^2 + |q|^2)p -(\a q +\frac{i}{2} q_x).
\eca.
\eeq 
Moreover, if $p$ is a solution of \eqref{ej}, then $\ti q= q-2i(\a-\bar \a) \frac{p}{|p|^2 -1}$ is a solution of the defocusing NLS. 
\een
\ethm 

Since $\li y(x,t), y(x,t)\ri_1$ may be zero at some points, BT of a smooth solution $q:\R^2\to \C$ of the defocusing NLS may develop singularities (for examples, cf. \cite{TU00}). 

B\"acklund transformations for the time-like $\ast$-MCF on $\R^{2,1}$ can be constructed in a similar manner as in Theorem \ref{cj}  and we get the following.

\bthm Let $\g$ be a solution of the time-like $\ast$-MCF, $\g=\frac{1}{2}\ast_\g H(\g)$, on $\R^{2,1}$, $g=(e_0, e_1, e_2)$, and $q$ the solution of the defocusing NLS as in Theorem \ref{fm}. Let $\a\in \C\bh \R$, and $y=(y_1,y_2)^t$ a solution of \eqref{ev}. Then 
$$\ti \g= \g+\frac{\Im(\a)}{|\a|^2||y||_1^2}\left(\frac{|y_1|^2+|y_2|^2}{2} e_0 -(y_1\bar y_2+ \bar y_1 y_2) e_1 +\frac{1}{i}(y_1\bar y_2-\bar y_1 y_2)e_2\right)$$
is a solution of the time-like $\ast$-MCF on $\R^{2,1}$, where $||y||_1^2= |y_1|^2- |y_2|^2$. 
\ethm

Permutability formula for BT and explicit soliton solutions for the time-like $\ast$-MCF can be obtained similarly as for the $\ast$-MCF on $\R^3$.  

\bs

\section{Space-like $\ast$-MCF on $\R^{2,1}$} \label{dwa}\ 

The outline of this section:
\ben
\item[(a)] We show that the differential invariants of solutions of the space-like $\ast$-MCF on $\R^{2,1}$ satisfy the second flow of the $SL(2,\R)$-hierarchy and use the B\"acklund transformations of the $SL(2,\R)$-hierarchy to construct B\"acklund transformations for the space-like $\ast$-MCF. 
\item[(b)] We note that the space-like normalized geometric Airy flow on $\R^{2,1}$ preserves the set 
$\calM_0$ of all space-like curves in $\R^{2,1}$ parametrized by arc-length and the principal curvature along one null parallel normal field is constant. Then the principal curvature along the other null parallel normal field satisfies the KdV.  So the normalized geometric Airy flow on $\calM_0$ gives a natural geometric interpretation of the KdV.   We use the BT for the KdV to construct BT for this restricted curve flow.
\een

For space-like solution of the $\ast$-MCF on $\R^{2,1}$, we identify $\R^{2,1}$ as $sl(2,\R)$ with
$$\li X, Y\ri = \frac{1}{2} \tr(XY).$$
Then 
$$\li X, X\ri = x_1^2 + x_2 x_3, \quad {\rm for\,\,} X=\bpm x_1 & x_2\\ x_3& -x_1\epm\in sl(2,\R).$$  Let 
\beq\label{}
a=\diag(1,-1), \quad b= \bpm 0& 1\\0&0\epm, \quad c= \bpm 0&0\\ 1&0\epm.
\eeq
Then $(a,b,c)$ is an ordered base of $sl(2,\R)$ satisfying
 \begin{align*}
 &\li a, a\ri=1, \quad \li b, c\ri=\frac{1}{2}, \quad \li b,b\ri=\li c, c\ri=0,\\
 & [a, b]= 2b, \quad [a, c]=-2c, \quad [b,c]= a.
 \end{align*}  We identify $sl(2,\R)$ as $\R^{2,1}$ by $xa+ yb+ zc\mapsto (x,y, z)^t$ with the bi-linear form  $x^2+ yz$.

Recall that the Hodge star operator on a two dimension Lorentzian space $(V, \li, \ri)$ is defined as follows:  Let $(e_1, e_2)$ be an ordered orthonormal basis of $V$, i.e., $\li e_1, e_2\ri =0$, $\li e_1, e_1\ri=-\li e_2, e_2\ri = 1$. The Hodge star operator $\ast:V\to V$ is defined by 
$$v_1\wedge(\ast v_2)= \li v_1, v_2\ri e_1\wedge e_2$$ for all $v_1, v_2\in V$.  It can be checked that 
$$\ast(e_1)= e_2, \quad \ast(e_2)=e_1.$$
So the Hodge star operator on the Lorentzian plane spanned by $b, c$ is
\beq\label{fn}
\ast b= b, \quad \ast c= -c.
\eeq
Since $\ast$ on $\C \simeq \R^2$ is the rotation by $\pi/2$ or the multiplication by $i=\sqrt{-1}$. So $\ast$ can be viewed as an analogue of the multiplication by $i$ in $\R^{1,1}$.  
  
Let 
\begin{align}
D_3&= \bpm 1& 0&0\\ 0&0 &\frac{1}{2}\\ 0 & \frac{1}{2} &0\epm, \label{cw}\\
O(3, D_3)&=\{A\in GL(3,\R)\n A^t D_3 A= D_3\}. \label{cw1}
\end{align} Then the Lie algebra of $O(3,D_3)$ is 
$$o(3, D_3)=\left\{\bpm 0 & -\frac{1}{2}r_2 & -\frac{1}{2} r_1\\ r_1 & r_3 &0\\ r_2 &0 & -r_3\epm\, \bigg|\, r_1, r_2, r_3\in \R\right\}.$$

If $\g:\R\to sl(2,\R)$ is  a space-like curve parametrized by the arc-length, then there is $g=(e_0, e_1, e_2):\R\to O(3, D_3)$ such that $e_0=\g_x$ and
$$g^{-1}g_x=\bpm 0& -\frac{k_2}{2}& -\frac{k_1}{2}\\ k_1 & 0 &0\\ k_2 & 0 &0 \epm.$$
So $(e_1, e_2)$ is a parallel null frame for the normal bundle $\nu(\g)$.  We call such $g$ a {\it null p-frame\/} for the space-like curve $\g$, and $k_i$ the principal curvature of $\g$ along the null parallel normal $e_i$ for $i=1, 2$. Note that the mean curvature vector is $H(\g)= (e_0)_x= k_1e_1 + k_2 e_2$. Again a proof similar to that of Theorem \ref{bb} gives 

\bthm\label{dp} Suppose $\g:\R^2\to  \R^{2,1}$ is a space-like solution of the $\ast$-MCF on $\R^{2,1}$, $\g=-\frac{1}{2}\ast_\g H(\g)$, parametrized by arc-length. Then there is $g=(e_0, e_1, e_2): \R^2\to O(3,D_3) $ such that
\ben
\item[(i)] $g(\cdot, t)$ is a null p-frame along $\g(\cdot, t)$ for all $t$, and
$$\bca g^{-1}g_x= \bpm 0& -\frac{k_2}{2} & -\frac{k_1}{2}\\ k_1 & 0& 0\\ k_2 &0&0\epm,\\
g^{-1}g_t= \frac{1}{4}\bpm 0& -k_2'& k_1'\\ -2k_1'& k_1k_2 &0\\ 2k_2' &0 & -k_1 k_2\epm,
\eca$$
\item[(ii)] $q= -\frac{k_1}{2}$ and $r= \frac{k_2}{2}$ satisfies the 2nd flow of the $SL(2,\R)$-hierarchy, i.e., 
\beq\label{bsa}
\bca q_t= -\frac{1}{2} (q_{xx} - 2q^2 r),\\
 r_t= \frac{1}{2} (r_{xx} - 2 qr^2).\eca
 \eeq 
\een
Moreover, if $\ti g:\R^2\to O(3, D_3)$ also satisfies (i) and (ii), then there is a constant $c\in \R\bh 0$ such that $\ti g= g\diag(1, c, c^{-1})$ and the corresponding null principal curvatures are related by $\ti k_1= c^{-1}k_1$ and $\ti k_2= ck_2$.  
\ethm

\brem If $(e_0, e_1, e_2)$ is a null p-frame along $\g$ and $c\in \R$ a non-zero constant, then $(e_0, ce_1, \frac{1}{c}e_2)$ is also a null p-frame and the corresponding null principal curvatures are related by $\ti k_1= c k_1$ and $\ti k_2= \frac{1}{c} k_2$.  In fact, if $(e_0, e_1, e_2)$ is a null p-frame along $\g$ then $(\g, e_0, e_1, e_2)$ is a $(O(1,1), V)$ moving frame along $\g$, where $O(1,1)= \{\diag(c, \frac{1}{c})\n c\in \R\bh 0\}$ and $V= e_{21}+\R (e_{32}- \frac{1}{2} e_{24}) \oplus \R( e_{42} -\frac{1}{2} e_{23})$.  Here we embed the Lie algebra of $\fR_1(3)$ as the following subalgebra of $gl(4,\R)$,
$$\left\{\bpm 0 & 0\\ y& A\epm\, \bigg|\, y\in \R^3, A\in o(3, D_3)\right\}.$$
\erem

\brem Theorem \ref{dp}  associate to each solution $\g$ of the space-like $\ast$-MCF on $\R^{2,1}$ an $O(1,1)$-orbit of solutions of \eqref{bsa}. The Pohlmeyer-Sym construction give the converse (the proof is the same as for the VFE). Let $\Psi_2$ denote the map from the orbit space of solutions of space-like $\ast$-MCF on $\R^{2,1}$ parametrized by arc-length under the rigid motion group $\fR_1(3)$ to the orbit space of solutions of \eqref{bsa} under the action of $O(1,1)$ defined by $\Psi_2([\g])= [u]$, where $u=qe_{12}+ re_{21}$ is a solution associated to $\g$ in Theorem \ref{dp}. Then $\Psi_2$ is a bijection and is the analogous Hasimoto transform for the space-like $\ast$-MCF on $\R^{2,1}$. 
\erem

Next we review the BT for the $SL(2,\R)$ hierarchy constructed in \cite{TU00}.  Let $\a_1, \a_2\in \R$,  $\{v_1, v_2\}$ a basis of  $\R^2$,  $\pi$ the projection of $\R^2$ onto $\R v_1$ along $\R v_2$, and $\pi'= \I-\pi$.  Let
$$h_{\a_1, \a_2, \pi}(\l)= \I +\frac{\a_1-\a_2}{\l-\a_1}\pi'.$$
A direct computation shows that 
$$h_{\a_1, \a_2,\pi}(\l)^{-1}= \I + \frac{\a_2-\a_1}{\l-\a_2} \pi'.$$

Note that the projection of $\pi$ of $\R^2$ onto $\R \bpm y_1\\ y_2\epm$ along $\R \bpm y_3\\ y_4\epm$ is 
$$\pi = \frac{1}{y_1y_4-y_2y_3}\bpm y_1y_4 & -y_1y_3\\ y_2y_4 & - y_2y_3\epm.$$

\bthm\label{fl} (\cite{TU00})  Let $u=q e_{12}+ re_{21}$ be a solution of the second flow \eqref{bsa} in the $SL(2,\R)$-hierarchy, and $\a\in \R$ a constant.  Then 
\beq\label{fp}
\bca \xi_x=-\bpm \a & q\\ r& -\a\epm \xi, \\
\xi_t= -\bpm \a^2-\frac{1}{2} qr & \a q -\frac{1}{2} q_x\\ \a r +\frac{1}{2} r_x & -\a^2 + \frac{1}{2} qr\epm \xi,
\eca
\eeq
is solvable.  Moreover, suppose  $\a_1\not=\a_2$ and $\xi_1=\bpm y_1\\ y_2\epm$ and $\xi_2=\bpm y_3\\ y_4\epm$ be solutions of \eqref{fp} with $\a= \a_1$ and $\a_2$ respectively. If $\xi_1, \xi_2$ are linearly independent, then 
$$\bca \ti q= q-\frac{2(\a_1-\a_2)}{\det(Y)} y_1 y_3, \\
\ti r= r- \frac{2(\a_1-\a_2)}{\det(Y)} y_2y_4,\eca$$ where $Y= (\xi_1, \xi_2)=\bpm y_1& y_3\\ y_2& y_4\epm$.
\ethm

\bthm Let $\g$ be a space-like solution of the $\ast$-MCF on $\R^{2,1}$, $\g_t= -\frac{1}{2} \ast_\g H(\g)$, parametrized by arc-length, $g=(e_0, e_1, e_2)$, and $q, r$ solution of the second flow \eqref{bsa} of the $SL(2,\R)$-hierarchy as in Theorem \ref{dp}.  Given $\a_1, \a_2\in \R$, let $Y=\bpm y_1& y_3\\ y_2& y_4\epm$ be as in Theorem \ref{fl}, then 
$$\ti \g= \g +\left(\frac{1}{\a_1}-\frac{1}{\a_2}\right)\frac{1}{ \det(Y)}\left(\frac{y_1y_4+ y_2y_3}{2} e_0 - y_1 y_3 e_1+ y_2 y_4 e_2  \right)$$
is a solution of the space-like $\ast$-MCF on $\R^3$.
\ethm

\ms
Next we discuss the space-like normalized geometric Airy flow on $\R^{2,1}$. 
A similar proof as for the geometric Airy flow on $\R^3$ gives the following.

\bthm If $\g$ is a solution of the space-like normalized geometric Airy flow on $\R^{2,1}$,
\beq\label{fr}
\g_t= \frac{1}{8} \li H(\g), H(\g)\ri e_0 +\frac{1}{4} \K_{e_0}^\perp H(\g),
\eeq
then there is $g=(e_0, e_1, e_2):\R^2\to O(3, D_3)$ such that $g(\cdot, t)$ is a null p-frame along $\g(\cdot, t)$ with principal curvatures $k_1,k_2$. Then $q=-\frac{k_1}{2}$ and  $r=\frac{1}{2} k_2$ satisfy the third flow of the $SL(2,\R)$-hierarchy, i.e.,  
\beq\label{eu}
\bca q_t= \frac{1}{4} (q_{xxx} - 6qrq_x),\\ r_t= \frac{1}{4}(r_{xxx} - 6qr r_x).\eca
\eeq
\ethm

It is known that \eqref{eu} leaves the condition $r=1$ invariant. Since $r=1$ means $k_2= 2$ in the above theorem, the space-like geometric Airy flow \eqref{fr} leaves the space curve with $k_2$ being a constant function invariant.
So we obtain the following.

\bthm\label{ew} Let $\calM_0$ denote the set of $\g:\R\to\R^{2,1}$ such that $||\g_x||=1$ and the principal curvature along  one  null parallel normal field is constant.  Then the space-like normalized geometric Airy flow \eqref{fr} on $\R^{2,1}$ leaves $\calM_0$ invariant. Moreover, if $\g$ is a solution of \eqref{fr} on $\calM_0$, then there is $g=(e_0, e_1, e_2):\R^2\to O(3,D_3)$ satisfying
\ben
\item[(i)] $g(\cdot, t)$ is a null p-frame for $\g(\cdot, t)$ and 
$$g^{-1}g_x= k_1(e_{21}-\frac{1}{2} e_{13}) + 2 (e_{31}-\frac{1}{2} e_{12}),$$
(so the null principal curvature $k_2$ along $e_2$ is equal to $2$),
\item[(ii)] $q=-\frac{1}{2} k_1$ satisfies the KdV, $q_t= \frac{1}{4}(q_{xxx} - 6 qq_x)$.
\een
\ethm

Hence the restriction of the normalized space-like geometric Airy flow \eqref{fr} on $\R^{2,1}$ to $\calM_0$ gives a natural geometric interpretation of the KdV. 

B\"acklund transformations for the KdV were constructed in several places. The version we need to use to construct BT for the curve flow \eqref{fr} on $\calM_0$ is proved in \cite{TU00}. 
Given $\xi, k\in \R$, let 
$$r_{\xi, k}(\l)= a\l + \bpm \xi & \xi^2-k^2\\ 1&\xi\epm, \quad a=\diag(1,-1).$$
Then we have
$$r_{\xi, k}(\l)^{-1}= \frac{ r_{-\xi, k}(\l)}{\l^2-k^2}.$$

Recall that the Lax pair for a solution $q$ of the KdV is 
$$\o_3= (a\l+ u)\rd x+ (a\l^3+ u\l^2+ Q_{-1}(u)\l+ Q_{-2}(u))\rd t,$$
where $a=\diag(1,-1)$, $ u= \bpm 0& q\\ 1&0\epm$, and 
$$Q_{-1}(u)= \frac{1}{2} \bpm -q & -q_x\\ 0& q\epm, \quad Q_{-2}(u)= \frac{1}{4}\bpm q_x& q_{xx}- 2q^2\\ -2q & -q_x\epm.$$
A frame $E(x,t,\l)$ of a solution $q$ of the KdV is a solution of $E^{-1}\rd E= \o_3$ and $E$ satisfies the {\it KdV reality condition},
\beq\label{}
\bca 
\overline{E(x,t,\bar\l)}=E(x,t,\l),\\ 
\phi(\l)^{-1}E(x,t,\l) \phi(\l) \, \, {\rm is\, an \, even\, function\, of\, \,}\l,
\eca
\eeq
where $\phi(\l)=\bpm 1& \l \\ 0& 1\epm$.

\bthm(\cite{TU00})
Let $E(x,t,\l)$ be a frame of the solution $q$ of the KdV, $c, \xi\in \R$, and 
$$\bpm y_1(x,t)\\ y_2(x,t)\epm= E(x,t, c)^{-1} \bpm c-\xi \\ 1\epm.$$
If $y_2\not=0$, let $\ti \xi = c-\frac{y_1}{y_2}$, then  
$$\ti q= -q + 2(\ti \xi^2 -c^2)$$
is a solution of the KdV. Moreover, 
$$\ti E(x,t,\l)=\frac{1}{\l^2-c^2}\, r_{\xi,c}(\l) E(x,t,\l)r_{-\ti\xi (x,t), c}(\l)$$
is a frame for the solution $\ti q$.  
\ethm

As a consequence we obtain BT for space-like normalized geometric Airy flow $\calM_0$.

\bthm Let  $\g$ be a solution of \eqref{fr} on $\calM_0$,  $g=(e_0, e_1, e_2)$ with principal curvatures $k_1$, $k_2=2$, and $q=-\frac{1}{2}k_1$ a solution of the KdV as in Theorem \ref{ew}.  Given a constant $c\in\R\bh 0$, then the following linear system for $y$ is solvable,
\beq\label{el}
\bca y_x=-\bpm c & q\\ 1& -c \epm y, \\
y_t= -\bpm c^3-\frac{c}{2}q +\frac{1}{4} q_x& c^2q -\frac{c}{2} q_x + \frac{1}{4}(q_{xx}- 2q^2)\\ c^2 - \frac{1}{2} q & -c^3+\frac{c}{2}q - \frac{1}{4} q_x   \epm\, y.\eca
\eeq
Moreover, if $y=(y_1, y_2)^t$ is a solution of \eqref{el}, then  
$$\ti \g= \g + \frac{1}{c^2} (-\ti \xi e_0 + (c^2-\ti \xi^2) e_1+ e_2)$$
is a solution of \eqref{fr} on $\calM_0$, where $\ti \xi= c-\frac{y_1}{y_2}$
\ethm

\bs
\section{Constructing integrable curve flows from Lax pairs}

In this section we explain how to use the Lax pair of a soliton equation to construct a geometric curve flow on a homogeneous space such that the evolution equation of its differential invariants is a soliton equation. We see that the $\ast$-MCF and the geometric Airy flow come up naturally this way. This method also gives the shape operator curve flow on Adjoint orbits, the Schr\"odingier flow on compact Hermitian symmetric spaces, and the central affine curve flows on $\R^n\bh 0$.  

We first review the {\it $G$-hierarchy\/} and its various restrictions (cf. \cite{DS84},  \cite{Sat84}, \cite{TU00}).

\ss\ni {\bf The $(G,a)$-hierarchy}\

Let $G$ be a complex simple Lie group, $\calG$ its Lie algebra, and  $a\in \calG$. Let $\calG_a^\perp$ denote the subspace perpendicular to the centralizer $\calG_a$ in $\calG$ with respect to the Killing form.  The {\it $(G,a)$-hierarchy\/} is the collection of flows on $C^\infty(\R, \calG_a^\perp)$ defined as follows:  It is known (cf., \cite{Sat84}, \cite{TU00}) that given a smooth $u:\R\to \calG_a^\perp$, there is a unique $Q(u,\l)\in C^\infty(\R, \calG)$ satisfying 
\beq\label{ey}
[\p_x+ a\l + u, Q(u,\l)]=0,\eeq
$Q(u,\l)$ is conjugate to $a\l$, and the power series expansion of $Q(u,\l)$ in $\l$ is
 $$Q(u,\l)=a\l+ u+ Q_{-1}(u)\l^{-1} + Q_{-2}(u)\l^{-2}+ \cdots.$$
 Moreover, entries of $Q_{-j}(u)$ are differential polynomial of $u$. Let $b\in \calG$ such that $b=\phi(a)$ for some analytic function $\phi$. {\it The flow in the $(G,a)$-hierarchy defined by $b\l^j$\/}  is 
 \beq\label{ez}
 u_t= [\p_x+ a\l + u, Q_{b, -(j-1)}(u)],
 \eeq
 where $Q_{b, -i}(u)$ is the coefficient of $\l^{-i}$ of the power series expansion of 
 $\phi(Q(u,\l)\l^{-1})\l^j$.  The Lax pair of \eqref{ez} is 
 \beq\label{fd}
 \o=(a\l+ u)\rd x + \left(b\l^j +\sum^{j-1}_{i=0} Q_{b, -i} \l^{j-1-i}(u)\right)\rd t.
 \eeq

\ss\ni {\bf The $(G,\sigma, a)$-hierarchy}\

Let $\sigma$ be an involution of $G$ such that the induced Lie algebra involution $\sigma_\ast$ on $\calG$ is complex linear, and $\calG_0$ and $\calG_1$ the eigenspace of $\sigma_\ast$ with eigenvalue $1, -1$ respectively. Then $\calG=\calG_0\oplus \calG_1$.  If $a\in \calG_1$, then the $(G,a)$-hierarchy leaves $C^\infty(\R, [a, \calG_1])$ invariant and $Q(u, \l)$ satisfies the $\sigma_\ast$ reality condition,
$$\sigma_\ast(Q(u, -\l))= Q(u,\l).$$
The restriction of the flows in the $(G,a)$-hierarchy to $C^\infty(\R, [a, \calG_1])$ is the {\it $(G,\sigma, a)$-hierarchy\/}. 
 
 \ss\ni {\bf The $(U,a)$-hierarchy}\
 
 Let $\tau$ be an involution of $G$ such that  the induced involution $\tau_\ast$ on the Lie algebra $\calG$ is conjugate linear.  
The fixed point set $U$ of $\tau$ in $G$ is a {\it real form of $G$\/}. 
  It is known that if $a\in \calU$, then the $(G,a)$-hierarchy leaves $C^\infty(\R, \calU_a^\perp)$ invariant, where $\calU_a^\perp$ is the orthogonal complement of the centralizer $\calU_a$ in $\calU$. Moreover, $Q_{b, -i}(u)\in \calU$.  The restriction of flows of the $(G,a)$-hierarchy to $C^\infty(\R, \calU_a^\perp)$ gives the {\it $(U,a)$-hierarchy\/}. 
 
 The Lax pair of flow \eqref{ez} for $u:\R^2\to \calU_a^\perp$ is given by the same formula \eqref{fd} and satisfies the {\it $U$-reality condition\/},
 $$\tau_\ast(\o(x,t,\bar\l))= \o(x,t,\l).$$
  A {\it frame\/} $E$ of a solution $u:\R^2\to \calU_a^\perp$ of \eqref{ez} is the solution of $E^{-1}\rd E= \o_\l$ satisfying the {\it $U$-reality condition\/}
 $$\tau(E(x,t,\bar\l))= E(x,t,\l).$$
 
 \ss\ni {\bf The $(\frac{U}{K}, a)$-hierarchy}\
 
 Let $U$ be the real form defined by the involution $\tau$, and $\sigma$ an involution of $G$ such that $\sigma_\ast$ is complex linear and  $\sigma\tau=\tau \sigma$. Let $K$ be the fixed point set of $\sigma$ in $U$, and  $\calU=\calK+\calP$ the Cartan decomposition (i.e., $\calP$ is the $-1$ eigenspace of $\sigma_\ast$ on $\calU$). If $a\in \calP$, then the flows in the $(\calU,a)$-hierarchy leaves $C^\infty(\R, [a,\calP])$ invariant and the resulting hierarchy is the {\it $(\frac{U}{K}, a)$-hierarchy}. The Lax pair \eqref{fd} for $u:\R^2\to [a, \calP]$ satisfies the {\it $\frac{U}{K}$-reality condition\/},
 $$\sigma_\ast(\o(x,t,-\l))= \o(x,t,\l), \quad \tau_\ast(\o(x,t,\bar\l))= \o(x,t,\l).$$
 A solution $E$ of $E^{-1}\rd E$ is a {\it frame\/} of a solution $u$ of the $(\frac{U}{K}, a)$-hierarchy if $E$ satisfies the {\it $\frac{U}{K}$-reality condition\/}
 $$\sigma(E(x,t,-\l))= E(x,t,\l)= \tau(E(x,t,\bar \l)).$$
 
\beg The $2\times 2$ AKNS, the $SU(2)$, the $SU(1,1)$, and the $SL(2,\R)$ hierarchies given in section \ref{bz} are the 
$$(SL(2,\C), a), \quad (SU(2),a), \quad (SU(1,1),a),  \quad (SL(2,\R),a_1)$$ hierarchies respectively, where $a=\diag(i, -i)$ and $a_1= \diag(1,-1)$.  
\eeg

\beg
 Suppose $a\in \calG$ is regular, i.e., $\calG_a=\calA$ is a Cartan subalgebra of $\calG$. The flow in the $(\calG, a)$-hierarchy generated by $b\l$ for $b\in \calA$ is the following $n$-wave equation defined by $G$ for $u:\R^2\to \calA^\perp$,
 \beq\label{fc}
 u_t= \ad(b)\ad(a)^{-1}(u_x)+ [u, \ad(b)\ad(a)^{-1}(u)]. 
 \eeq
 \eeg
 
 \ss\ni {\bf Integrable curve flows on Adjoint orbits}\
 
 We give a general method of using the Lax pair of a soliton equation to construct curve flows on the Adjoint $U$-orbit $M$ at $a$ in $\calU$ whose differential invariants satisfies the given soliton equation. For the  $(U,a)$-hierarchy, we have the following (cf. [TU06], [TerTho01])
  \ben
 \item[(i)] Given a solution $u$ of \eqref{ez} of the $(U,a)$-hierarchy, let $g:\R^2\to U$ be a solution of 
 \beq\label{fa}
 \bca g^{-1}g_x= u,\\ g^{-1}g_t= Q_{b, -(j-1)}(u),\eca
 \eeq
 i.e., $g$ is a parallel frame for the flat connection $\o(\cdot, \cdot, 0)$, the Lax pair \eqref{fd} at $\l=0$.
  Set 
  $$\g(x,t)= g(x,t) ag^{-1}(x,t).$$ Then $\g(\cdot, t)$ is a curve on $M=\Ad(U)(a)$ for all $t$ and $g(\cdot, t)$ is an Adjoint moving frame along $\g(\cdot, t)$, and $u=g^{-1}g_x\in \calU_a^\perp$ is the differential invariants (see Example \ref{fs}).  
 \item[(ii)] The second equation of \eqref{fa} implies that  
 \beq\label{fb}
 \g_t= g[Q_{b,-(j-1)}(u), a]g^{-1}.\eeq
 Note that this curve flow does not preserve the arc-length and has no special parameters. So the number of differential invariants defined by the Adjoint moving frame should equal to $\dim(M)$.
 \item[(iii)] Since entries of $Q_{b, -(j-1)}(u)$ are differential polynomials in $u$, the right hand side of \eqref{fb} is a quantity written purely in terms of the moving frame and the differential invariants of a curve.  We use the same proof as for the $SU(2)$-hierarchy to see that 
 $$Q(h u h^{-1})= hQ(u)h^{-1}$$ for a constant $h\in U_a$ and $u\in C^\infty(\R, \calU_a^\perp)$. So $U_a$ acts on the space of solutions of \eqref{ez} of the $(U,a)$-hierarchy. It can be checked that the right hand side of \eqref{fb} is independent of the choice of Adjoint moving frames.  So \eqref{fb} is a geometric curve flow on the Adjoint orbit $M$.  
 \item[(iv)] Conversely, given a solution $\g:\R^2\to M$ of the curve flow \eqref{fb} on $M$, it was proved in \cite{TU06} and \cite{TerTho01} that there is $g:\R^2\to U$ such that $g(\cdot, t)$ is an Adjoint frame of $\g(\cdot, t)$ and $u:=g^{-1}g_x$ is a solution of the \eqref{ez} of the $(U,a)$-hierarchy. Moreover, $[\g]\mapsto [u]$ is a bijection from the orbit space of solutions of the curve flow \eqref{fb} under the action of $U$ and the orbit space of solutions of \eqref{ez} under the action of $U_a$. 
  \een
This construction gives natural integrable curve flows on Adjoint $U$-orbits.  
  For example, 
 \ben
 \item[$\di$]
 if $a\in \calU$ is regular, then the flow \eqref{ez} in the $(U,a)$-hierarchy defined by $b\l$ is the $n$-wave equation \eqref{fc} and the curve flow \eqref{fb} is the shape operator curve flow $\g_t= A_{\hat b}(\g_x)$ on principal Adjoint orbit $M=\Ad(U)(a)$ in $\calU$ (cf. \cite{Fer95}, \cite{TerTho01}), where $\hat b$ is the normal field on $M$ defined by $\hat b(gag^{-1})= gbg^{-1}$. 
 \item[$\di$] if $U$ is compact and $a\in U$ such that $\ad(a)^2=-{\rm id}$ on $\calU_a^\perp$, then $M=\Ad(U)(a)$ is a Hermitian symmetric space. The flow in the $(U,a)$-hierarchy generated by $a\l^2$ is the $\frac{U}{U_a}$-NLS hierarchy constructed in \cite{ForKul83}. It was proved in \cite{TU06} that the curve flow \eqref{fb} for $a\l^2$ on $M$ is the Schr\"dingier flow on the Hermitian symmetric space $M$.
 \een
 
 \ms\ni {\bf Integrable curve flow on $\calU$ with constraint}\
 
 We explain below a method for constructing integrable curve flows on the flat space $\calU$ with constraint from the soliton equation \eqref{ez} of the $(U,a)$-hierarchy, where $U$ is a compact Lie group, $\li a, a\ri=1$, and $\li\, ,\ri$ is an ad-invariant inner product on $\calU$. 
\ben
\item[(a)] Let $E$ be a frame of a solution $u$ of \eqref{ez} in the $(U,a)$-hierarchy defined by $b\l^j$, and $\g= \frac{\p E}{\p\l}E^{-1}|_{\l=0}$. A direct computation as for Theorem \ref{bda} shows that
\begin{align}
&\g_x= \phi a\phi^{-1}, \label{ff}\\
& \g_t= \phi Q_{b, -(j-2)}(u)\phi^{-1}, \label{fg}
\end{align}
where $\phi(x,t)= E(x,t,0)$. So
\beq\label{fz}
\phi^{-1}\phi_x= u, \quad \phi^{-1}\phi_t= Q_{b, -(j-1)}(u),
\eeq
where $u\in \calU_a^\perp$. Let $a_0= a, \ldots, a_m$ be an orthonormal basis of $\calU_a^\perp$ with respect to $\li\, ,\ri$ of $\calU$. Let $e_i= \phi a_i \phi^{-1}$ for $0\leq i\leq m$. Then $g=(e_0, \ldots, e_m)$ is an orthonormal frame along $\g$ with $e_0=\g_x$. It follows from \eqref{fz} that we have $$(e_i)_x= \phi [u, a_i] \phi^{-1}$$ for all $0\leq i\leq m$.  So we can use Lie theory to write down the metric invariants of the curve in terms of $u=\phi^{-1}\phi_x$. 
\item[(b)] Write $\phi Q_{b, -(j-2)}(u)\phi^{-1}$ in terms of $g$ and the metric invariants and show that it is a geometric quantity (i.e., independent of the choice of frames). Note that the curve flow \eqref{fg} is only defined on 
\beq\label{fe}
\calM_a=\{\g:\R\to \calU\n ||\g_x||=1, \g_x\in \Ad(U)(a)\}.
\eeq
The constraint condition $\g_x\in \Ad(U)(a)$ implies that we only expect the number of local invariants of $\g$ should be $(m+1)$ functions where $m+1=\dim(\calU_a^\perp)$. Note that the metric invariants are given by $u$ which is a $\calU_a^\perp$ valued function, so the number of  local invariants for $\g\in \calM_a$ is correct. 
\item[(c)] The map $\Phi$ from the orbit space of  of \eqref{ez} under the action of $U_a$ to  the orbit space of solutions the curve flow \eqref{fg} under the action of $U$ defined by $\Phi([u])=[\g]$ is well-defined. 
\item[(d)] If the map $\Phi$ is a bijection, then we can use soliton theory of the $(U,a)$-hierarchy to study the curve flow \eqref{fg} on $\calM_a$. 
  \een
   
 The following are some known results.
  \ben
 \item If $\g$ is a solution of the curve flow \eqref{fg} on $\calM_a$ then $\a= \g_x$ is a solution of \eqref{fb} on the Adjoint orbit $\Ad(U)(a)$. 
  \item The methods described above can be applied  to the $(U,a)$-hierarchy and the $(\frac{U}{K},a)$-hierarchy as well. 
 \item  We have carried out the above program for curves in $\calU$ in this paper when the rank of $U$ is one. In this case, the class of curves $\calM_a$ is the space of all curves parametrized by arc-length.  So the curve flow \eqref{fg} is a flow on all immersed curves on $\calU$.  In fact, this is how the $\ast$-MCF and the geometric Airy flow are discovered. 
 \item Thorbergsson and the author use this method in \cite{TerTho14} to show that the third curve flow \eqref{fg} constructed from $(S^n, a)$-hierarchy is the normalized geometric Airy flow on $\R^n$. New integrable curve flows on $\C^n$ and $\H^n$ whose differential invariants satisfy the flows in the $(\C P^n,a)$- and the $(\H P^n, a)$- hierarchies are also studied in \cite{TerTho14}.
 \item When $\frac{U}{K}$ is a Hermitian symmetric space, the curve flow \eqref{fg} on $\calM_a$ was studied in \cite{LanPer00}. 
 \een  
 
 \ss\ni {\bf Central affine curve flows on $\R^n\bh 0$}\
 
Drinfeld-Sokolov  constructed in \cite{DS84} a KdV type hierarchy for each affine Kac-Moody algebra and showed that the one associated to $\hat A_1^{(1)}$ is the KdV hierarchy and the one associated to $\hat A_{n-1}^{(1)}$ for $n>2$ is the Gelfand-Dickey hierarchy.  The Lax pair for KdV used in \cite{DS84} is 
\beq\label{fh}
\Theta(\cdots, z)=\bpm 0& q+z\\ 1&0\epm \rd x + \bpm \frac{1}{4}q_x & z^2 +\frac{q}{2} z +\frac{1}{4} (q_{xx}- 2q^2)\\ z-\frac{1}{2}q & -\frac{1}{4}q_x\epm\rd t.
\eeq
This is gauge equivalent to the Lax pair $\o$ given by \eqref{fq}. In fact, 
$$\phi(\l) \o(x,t,\l)\phi(\l)^{-1}= \Theta(x,t, \l^2),$$
where $\phi=\bpm 1& \l \\ 0&1\epm$.  Let $q$ be a solution of the KdV, and $g:\R^2\to SL(2,\R)$ a parallel frame for $\Theta(\cdots, 0)$, i.e.,
\beq\label{fi}
\bca g^{-1}g_x= \bpm 0&q\\ 1&0\epm,\\ g^{-1}g_t= \bpm \frac{q_x}{4} & \frac{1}{4}(q_{xx}-2q^2)\\ -\frac{1}{2} q& -\frac{1}{4}q_x\epm.\eca
\eeq
Let $\g$ denote the first column of $g$. Then the first equation of \eqref{fi} implies that $g=(\g, \g_x)$ and $\g_{xx}= q\g$. It follows from Example \ref{ft} that  $\g(\cdot, t)$ is a curve on $\R^2\bh 0$ parametrized by the central affine arc-length parameter, $g(\cdot, t)$ is the central affine moving frame along $\g(\cdot, t)$, and $q(\cdot, t)$ is the central affine curvature of $\g(\cdot, t)$ for all $t$.  The second equation of \eqref{fi} implies that 
\beq\label{fj}
\g_t= \frac{1}{4}q_x \g -\frac{1}{2}q \g_x.\eeq
This is the curve flow on $\R^2\bh 0$ considered by Pinkall in \cite{UP95}.  It can be checked that $q\mapsto \g$ is a bijection from the space of solutions of the KdV to the space of solutions of the curve flow \eqref{fj} on $\R^2\bh 0$ parametrized by the central affine arc-length parameter.  There have been several works on the curve flow \eqref{fj} on $\R^2\bh 0$ concerning the higher flows and the bi-Hamiltonian structure (cf. \cite{CIM09}, \cite{FK10}, \cite{FK13}, \cite{TW14a}). Wu and the author also constructed B\"acklund transformations  and proved  the periodic Cauchy problem for \eqref{fj} has long time existence in \cite{TW14a}.   

Note that the two gauge equivalent Lax pairs of the KdV give rise to two different geometric interpretations of the KdV: One as space-like geometric Airy flow on $\R^{2,1}$ with constant principal curvature along a null parallel normal field and the other as the central affine curve flow on $\R^2\bh 0$.  

The Lax pair of the second flow in the $A^{(1)}_{n-1}$-KdV hierarchy for $n>2$ is of the form,
$$\o= (a\l + b+ u) \rd x+ ((e_{1, n-1}+ e_{2,n})\l + Z_{2,0}(u)) \rd t,$$
where $b=\sum_{i=1}^{n-1} b_{i+1, i}$ and $u=\sum_{i=1}^{n-1} u_i e_{in}$. Let $g$ be a solution of 
\beq\label{ga} g^{-1}g_x= b+u,\quad g^{-1}g_t= Z_{2,0}(u),\eeq
and $\g$ the first column of $g$. The first equation of  \eqref{ga} implies that $\g\in \calM_n(\R)$, $g=(\g, \g_x, \ldots, \g_x^{(n-1)})$ is the central affine moving frame along $\g$, and $u_1, \ldots, u_{n-1}$ are the central affine curvatures, where $\calM_n(\R)$ is defined in Example \ref{ft}. We use the second equation of \eqref{ga} to write down the following curve flow  on $\R^n\bh 0$ and show that its central affine curvatures satisfies the 2nd flow in the $A^{(1)}_{n-1}$-KdV hierarchy,
 \beq\label{ncf2}
\gamma_t=-\frac{2}{n}u_{n-1}\gamma+\gamma_{xx}.
\eeq
Higher order commuting central affine curve flows and a Poisson structure for \eqref{ncf2} on $\R^3\bh 0$ were given in \cite{CIM13}.  A bi-Hamiltonian structure, higher order curve flows,  and global periodic Cauchy problem for \eqref{ncf2} on $\R^n\bh 0$ are given in \cite{TW14b}.  B\"acklund transformations for \eqref{ncf2} on $\R^n\bh 0$ are constructed in \cite{TW14c}. 

\ms 

Since the $\ast$-MCF and the geometric Airy flow on $\R^3$ and $\R^{2,1}$ are integrable, the following questions arise naturally:
\ben
\item[$\bu$]
Does the periodic Cauchy problem of the  the $\ast$-MCF on $(N^3,g)$ have long time existence?
\item[$\bu$] Does the periodic Cauchy problem of the geometric Airy flow on a complete Riemannian manifold have long time existence?
\item[$\bu$] Is the $\ast$-MCF on a $3$ dimension homogeneous Riemannian $G$-manifold integrable?
\item[$\bu$] Let $(N,\rg)$ be a Riemannian manifold, $\dim(M)= \dim(N)-2$, and $\calI(M, N)$ the space of immersions from $M$ to $N$.  The $\ast$-MCF  on $\calI(M, N)$ is the following equation $$\frac{\p f}{\p t}= \ast_f(H(f(\cdot, t)))$$ for $f:\R\times M\to N$, where $f(\cdot, t)\in \calI(M,N)$, $H(f(\cdot, t))$ is the mean curvature vector field of the immersion $f(\cdot, t)$, and $\ast_f$ is the Hodge star operator on the normal plane $\nu(M)_f$.  Does the Cauchy problem of the $\ast$-MCF on $\calI(M,N)$ have long time existence?
\een

\bs


\begin{thebibliography}{99}

\bibitem{AblCla91}
Ablowitz, M.J.,Clarkson, P.A.,\emph{{S}olitons,
non-linear evolution equations and inverse scattering},  Cambridge Univ.
Press (1991)

\bibitem{BeaCoi84}
Beals, R., Coifman, R.R.,\emph{{S}cattering and inverse scattering for
first order systems}, Commun. Pure Appl. Math. \textbf{37} (1984), 39--90

\bibitem{BeaCoi85} 
Beals, R., Coifman, R.R., \emph{{I}nverse scattering and evolution
equations}, Commun. Pure Appl. Math., \textbf{38} (1985), 29-42

\bibitem{Bour93}
Bourgain, J., \emph{Fourier transform restriction phenomena for certain lattice subsets and applications to nonlinear evolution equations. I. Schršdinger equations},
Geom. Funct. Anal. 3 (1993), 107Ð156. 

\bibitem{CaIv05}
Calini, A., Ivey, T., \emph{Finite-Gap Solutions of the Vortex Filament Equation: Genus One Solutions and Symmetric Solutions}, J. Nonlinear Sci.
\textbf{15} (2005) pp. 321Ð 361.

\bibitem{CaIv07}
Calini, A., Ivey, T., \emph{Finite-gap solutions of the vortex filament equation: isoperiodic deformations}, J. Nonlinear Sci. \textbf{17} (2007) 527Ð567. 

\bibitem{CIM09}
Calini, A., Ivey, T., Mar{\'{\i}} Beffa,G., \emph{Remarks on KdV-type flows on star-shaped curves}, Phys. D 238 (2009), 788Ð797.

\bibitem{CIM13}
Calini, A., Ivey, T., Mar{\'{\i}} Beffa,G., \emph{Integrable flows for starlike curves in centroaffine space}, SIGMA Symmetry Integrability Geom. Methods Appl. \textbf{9} (2013), 022, 21pp.

\bibitem{DoSa94}
 Doliwa, A., Santini, P. M. , \emph{An elementary geometric characterization of the integrable motions of a curve}, Phys. Lett. A  \textbf{185} (1994)  373Ð384.
 
\bibitem{DS84}
Drinfel'd, V.G., Sokolov, V.V., \emph{Lie algebras and equations of Korteweg-de Vries type},  (Russian) Current problems in mathematics, \textbf{24} (1984),  81--180, Itogi Nauki i Tekhniki, Akad. Nauk SSSR, Vsesoyuz. Inst. Nauchn. i Tekhn. Inform., Moscow.

\bibitem{FO99}
 Fels, M., Olver, P.J., \emph{Moving coframes, II: Regularization and theoretical foundations}, Acta Appl. Math. 55:2 (1999), 127Ð208.

\bibitem{Fer95} Ferapontov, E.V., \emph{Isoparametric hypersurfaces in spheres, integrable non-diagonalizable systems of hydrodynamic type, and N-wave systems}, Differential Geom. Appl. 5:4 (1995), 335Ð369

\bibitem{ForKul83}
Fordy, A.P., Kulish, P.P., \emph{{N}onlinear Schr\"odinger equations and simple
Lie algebra}, Commun. Math. Phys., \textbf{89} (1983), 427-443

\bibitem{FK10}
Fujioka, A., Kurose, T., \emph{{H}amiltonian formalism for the higher KdV flows on the space of closed complex equicentroaffine curves}, Int. J. Geom. Methods Mod. Phys. \textbf{7(1)} (2010), 165-175.

\bibitem{FK13}Fujioka, A., Kurose, T., \emph{{M}ulti-Hamiltonian structures on space of closed equi-centroaffine plane curves associated to higher KdV flows}, preprint,  Arxiv: math.dg 1310.1688

\bibitem{Has72}
Hasimoto, H., \emph{ A soliton on a vortex filament}, J. Fluid Mech. 51 (1972), 477Ð485.

\bibitem{Its76}
Its, A. R., \emph{Inversion of hyperelliptic integrals, and integration of nonlinear differential
equations}, (Russian. English summary) Vestnik Leningrad. Univ. 1976, no. 7 Mat. Meh. Astronom. vyp. 2, 39Ð46, 162.

\bibitem{LanPer91}
Langer, J., Perline, R., \emph{ Poisson geometry of the filament equation}, J. Nonlinear Sci. \textbf{1} (1991) 71Ð93.

\bibitem{LanPer00}
Langer, J., Perline, R., \emph{Geometric realizations of Fordy-Kulish non- linear Schršdinger systems}, Pacific J. Math. \textbf{195} (2000), 157Ð178. 

\bibitem{MB08a} Mar{\'{\i}} Beffa,G., \emph{Geometric realizations of bi-Hamiltonian completely integrable systems},  SIGMA Symmetry Integrability Geom. Methods Appl. 4 (2008), Paper 034, 23 pp.

\bibitem{MB08b} Mar{\'{\i}} Beffa,G., \emph{Projective-type differential invariants and geometric curve evolutions of KdV-type in flat homogeneous manifolds},  Ann. Inst. Fourier (Grenoble) 58 (2008) 1295Ð1335.

\bibitem{MB09a} Mar{\'{\i}} Beffa,G., \emph{Hamiltonian evolution of curves in classical affine geometries}, Phys. D 238 (2009) 100Ð115.

\bibitem{MB10a} Mar{\'{\i}} Beffa,G., \emph{Bi-Hamiltonian flows and their realizations as curves in real semisimple homogeneous manifolds},  Pacific J. Math. 247 (2010), 163Ð188.

\bibitem{Pal60}
Palais, R.S., \emph{Slices and equivariant embeddings}, a chapter in Seminar on transformation groups by A. Borel,  Annals of Mathematics Studies,\textbf{46} (1960), 101-115.


\bibitem{Poh76}
Pohlmeyer, K., \emph{Integrable Hamiltonian systems and interactions through quadratic constraints}, Comm. Math. Phys. 46 (1976) p 207.
 
\bibitem{UP95}
Pinkall, U., \emph{Hamiltonian flows on the space of star-shaped
curves},  Results Math. \textbf{27(3-4)} (1995), 328--332. 

\bibitem{SaYa98}
Sasaki, N., Yasui, Y., \emph{ Differential geometry of the vortex filament equation}, J. Geom. Phys. \textbf{28} (1998)  p 195

\bibitem{Sat84}
Sattinger, D.H., \emph{{H}amiltonian hierarchies on semi-simple Lie
algebras}, Stud. Appl. Math., \textbf{72} (1984), 65--86

\bibitem{Sym85} 
Sym, A., \emph{ Soliton surfaces and their applications, in: Geomet- rical aspects of the Einstein equations and integrable systems}, Lecture Notes in Physics \textbf{239} (1985) p 154.

\bibitem{Ter97}
Terng, C.L., \emph{{S}oliton equations and differential
geometry}, J. Differential Geometry, \textbf{45} (1997), 407--445

\bibitem{TerTho01} 
Terng, C.L., Thorbergsson, G., \emph{Completely integrable curve flows on adjoint orbits}, Results Math. 40 (2001), 286Ð309.  Dedicated to Shing-Shen Chern on his 90th birthday. 

\bibitem{TerTho14} 
Terng, C.L., Thorbergsson, G., \emph{Integrable curve flows on $\C^n$ and $\H^n$}, preprint.

\bibitem{TerUhl98}
Terng, C.L., Uhlenbeck, K., \emph{{P}oisson actions and
scattering theory for integrable systems}, Surveys in Differential
Geometry: Integrable systems (A supplement to J. Differential
Geometry), \textbf{4} (1998), 315--402

\bibitem{TU00}
Terng, C.L., Uhlenbeck, K., \emph{B\"acklund transformations and
loop group actions}, Comm. Pure Appl. Math. \textbf{53} (2000), 1--75.

\bibitem{Ter03}
Terng, C.L., \emph{{G}eometries and symmetries of soliton equations and integrable elliptic systems}, to appear in 
Surveys on Geometry and Integrable Systems,  Advanced Studies in Pure Mathematics, Mathematical Society
of Japan, math.DG/0212372

\bibitem{TU06}
Terng, C.L., Uhlenbeck, K., \emph{Schršdinger flows on Grassmannian- ansÓ, pp. 235Ð256 in Integrable systems}, geometry, and topology, edited by C.-L. Terng, AMS/IP Stud. Adv. Math. \textbf{36} (2006), Amer. Math. Soc., Providence, RI.

\bibitem{TerWan08}
Terng, C.L., Wang, E., \emph{Transformations of flat Lagrangian immersions and Egoroff nets},  Asian J. Math. \textbf{12}  (2008) 99Ð119

\bibitem{TW14a}
Terng, C.L., Wu, Z., \emph{Central affine curve flow on the plane}, to appear in JFPTA Festschrift Volume for Mme Choquet-Bruhat (2014), 22pp.

\bibitem{TW14b}
Terng, C.L., Wu, Z., \emph{n-dimensional central affine curve flow}, preprint.

\bibitem{TW14c}
Terng, C.L., Wu, Z., \emph{B\"acklund transformations for the n-dim central affine curve flow}, preprint.

\bibitem{Uhl89}
Uhlenbeck, K., \emph{{H}armonic maps into Lie group (classical solutions of the
Chiral model)}, J. Differential Geometry, \textbf{30} (1989), 1--50

\bibitem{ZS72}
Zakharov, V., Shabat, A., \emph{Exact theory of two-dimensional self-focusing 
and one-dimensional self-modulation of waves in nonlinear media}, Soviet Physics JETP 34:1, 62-69 (1972). 


\end{thebibliography}
\end{document}